\documentclass[12pt]{amsart}
\usepackage{amssymb}
\input{diagrams.tex}
\diagramstyle[scriptlabels,height=8mm,width=8mm]
\theoremstyle{plain}

\newtheorem{thm}{Theorem}[section]
\newtheorem{cor}[thm]{Corollary}
\newtheorem{lem}[thm]{Lemma}
\newtheorem{prop}[thm]{Proposition}

\theoremstyle{definition}
\newtheorem{defi}[thm]{Definition}
\newtheorem{conj}[thm]{Conjecture}
\newtheorem{conv}[thm]{Convention}
\newtheorem{nota}[thm]{Notation}
\newtheorem{rem}[thm]{Remark}
\newtheorem{rems}[thm]{Remarks}
\newtheorem{exa}[thm]{Example}
\newtheorem{exas}[thm]{Examples}
\newtheorem{sit}[thm]{}

\newcommand{\brem}{\begin{rem}}
\newcommand{\brems}{\begin{rems}}
\newcommand{\erem}{\end{rem}}
\newcommand{\erems}{\end{rems}}
\newcommand{\bexa}{\begin{exa}}
\newcommand{\bexas}{\begin{exas}}
\newcommand{\eexa}{\end{exa}}
\newcommand{\eexas}{\end{exas}}
\newcommand{\bdefi}{\begin{defi}}
\newcommand{\edefi}{\end{defi}}
\newcommand{\bcor}{\begin{cor}}
\newcommand{\ecor}{\end{cor}}
\newcommand{\blem}{\begin{lem}}
\newcommand{\elem}{\end{lem}}
\newcommand{\bconv}{\begin{conv}}
\newcommand{\econv}{\end{conv}}
\newcommand{\bconj}{\begin{conj}}
\newcommand{\econj}{\end{conj}}
\newcommand{\bprop}{\begin{prop}}
\newcommand{\eprop}{\end{prop}}
\newcommand{\bthm}{\begin{thm}}
\newcommand{\ethm}{\end{thm}}
\newcommand{\bnota}{\begin{nota}}
\newcommand{\enota}{\end{nota}}
\newcommand{\bsit}{\begin{sit}}
\newcommand{\esit}{\end{sit}}
\newcommand{\be}{\begin{eqnarray}}
\newcommand{\ee}{\end{eqnarray}}
\newcommand{\bproof}{\begin{proof}}
\newcommand{\eproof}{\end{proof}}
\newcommand{\no}{\noindent}
\def\ba{\begin{array}}
\def\ea{\end{array}}

\def\lto{\longrightarrow}
\def\hto{\hookrightarrow}

\def\fm{{\mathfrak m}}

\def\cO{{\mathcal O}}

\newcommand{\trdeg}{\operatorname{trdeg}}
\newcommand{\Sing}{\operatorname{Sing}}
\newcommand{\Spec}{\operatorname{Spec}}
\newcommand{\Frac}{\operatorname{Frac}}

\newcommand{\Proj}{\operatorname{Proj}}

\newcommand{\id}{\operatorname{id}}
\newcommand{\Der}{\operatorname{Der}}

\newcommand{\Aut}{{\operatorname{Aut}}}

\newcommand{\SL}{{\bf {SL}}}
\newcommand{\GL}{{\bf {GL}}}
\newcommand{\PSL}{{\bf {PSL}}}
\def\MMh{{\rm MM_h}}
\def\MM{{\rm MM}}
\def\PGL{{\bf PGL}}
\def\fg{{\mathfrak g}}
\renewcommand{\div}{{\operatorname{div}}}

\def\sl{{\mathfrak{sl}}}

\newcommand{\A}{{\mathbb A}}
\newcommand{\pP}{{\mathbb P}}
\newcommand{\R}{{\mathbb R}}
\newcommand{\C}{{\mathbb C}}
\newcommand{\Q}{{\mathbb Q}}
\newcommand{\Z}{{\mathbb Z}}
\newcommand{\N}{{\mathbb N}}
\newcommand{\T}{{\mathbb T}}

\newcommand{\p}{{\partial}}

\def\bO{{\bar O}}
\def\quot{/\hskip-3pt/}

\title[Locally nilpotent derivations on $\C^*$-surfaces]{Locally
nilpotent derivations on affine surfaces with
a $\C^*$-action}

\author{Hubert Flenner}
\address{Fakult\"at f\"ur Mathematik,
Ruhr Universit\"at Bochum,
Geb.\ NA 2/72,
Universit\"ats\-str.\ 150,
44780 Bochum, Germany}
\email{Hubert.Flenner@ruhr-uni-bochum.de}

\author{Mikhail Zaidenberg}
\address{Universit\'e
Grenoble I, Institut Fourier, UMR 5582 CNRS-UJF, BP 74,
38402 St.\ Martin
d'H\`eres c\'edex, France}
\email{zaidenbe@ujf-grenoble.fr}

\thanks{
{\bf Acknowledgements:} This research was done during a visit of
the first author at the Fourier Institute of the University of
Grenoble, a stay of both authors at the Max
Planck Institute of Mathematics at Bonn and of the second author
at the Ruhr University at Bochum. The authors thank these
institutions for their generous support and excellent working
conditions.}

\thanks{
\mbox{\hspace{11pt}}{\it 1991 Mathematics Subject
Classification}:
14R05, 14R20, 14J50.\\
\mbox{\hspace{11pt}}{\it Key words}: $\C^*$-action,
$\C_+$-action,
graded algebra, affine surface, cyclic quotient
singularity}

\begin{document}

\begin{abstract}
We give a classification of normal affine surfaces admitting an
algebraic group action with an open orbit. In particular an
explicit algebraic description of the affine coordinate rings and
the defining equations of such varieties is given. By our methods
we recover many known results, e.g.\ the classification of normal
affine surfaces with a `big' open orbit of Gizatullin \cite{Gi1,
Gi2} and Popov \cite{Po} or some of the classification results of
Danilov-Gizatullin \cite{DaGi2}, Bertin \cite{Be1,Be2} and others.

\end{abstract}

\maketitle

\centerline{\it Dedicated to
M.\ Miyanishi on occasion of his
60th birthday}

{\footnotesize \tableofcontents}

\section*{Introduction}

Let $G$ be an algebraic group acting on a normal affine algebraic
surface $V$. By  classical results of Gizatullin \cite{Gi1} and
Popov \cite{Po}, if $V$ is smooth and $G$ has a big open orbit
$O$ (that is, $V\backslash O$ is finite), then $V$ is one of the
surfaces
$$
\C^{*2},\quad \A_\C^2, \quad \C^*\times\A_\C^1,
\quad \pP^1\times \pP^1\backslash
\Delta,
\quad \pP^2\backslash\bar\Delta,
$$
where $\Delta\subseteq \pP^1\times\pP^1$ is the diagonal and
$\bar\Delta\subseteq\pP^2$ is a nondegenerate quadric.
Furthermore, if $V$ is singular then $V\cong V_d$ is  the  Veronese
cone  $\A_\C^2/\Z_d$, where $\Z_d$ acts on $\A_\C^2$ via
multiplication with the group of $d$-th roots of unity (see
Example \ref{pop}).

The aim of this paper is to give more generally {\em a
description of all normal affine surfaces $V=\Spec A$ (over the
ground field $\C$) that admit an action of an algebraic group
with an open orbit}.  As was suggested by Popov \cite{Po} and
confirmed in the smooth case by Bertin \cite{Be2}, either such
a surface $V$ is isomorphic to $\C^{*2}$, or  a semidirect
product of $\C^*$ and $\C_+$ acts on $V$ with an open orbit
(Proposition
\ref{agoo}). We provide a classification of all  such surfaces
  in Section 3. This leads to a new proof of the
Gizatullin-Popov Theorem above (see Section 4.4) which uses
only elementary facts from Lie theory. For generalizations of
this result see also \cite{Ak2, HO}.

Our interest in such actions is inspired by the role that they
play in certain classification problems, e.g.\  in the proof
of linearization of regular $\C^*$-actions on $\A^3_\C$
\cite{KKMLR}. Usually in applications, to an affine variety $V$
with a $\C_+$-action one associates (non canonically) another
one, say, $V'$ with a $\C^*$- and $\C_+$-action (see e.g.,
\cite{ML1, Za} and Remark \ref{adlnd}.3 below). Therefore it is
of particular importance to classify such varieties.

$\C^*$-actions on algebraic surfaces were extensively studied in
the literature, see \cite{FlZa1} and the references given
therein, and also \cite{AlHa} for a  generalization to higher
dimensions. A $\C^*$-action on $V$  gives rise to a grading
$A=\bigoplus_{i\in\Z}A_i$.  We will rely here on our previous
paper \cite{FlZa1}  to describe the graded components $A_i$  in terms
of the Dolgachev-Pinkham-Demazure construction (the {\it
DPD construction}, for short).


Classification results for $\C_+$-actions on affine surfaces can
be found in \cite{Dani, Fi, ML2, Wi, MaMiy},
\cite{BaML2}-\cite{BaML3}, \cite{DaiRu} and \cite{Du1, Du2}. It
is well known \cite{Ren} that a $\C_+$-action gives rise to a
locally nilpotent derivation $\p$ of $A$ (see Proposition
\ref{Ren}). The condition that a semidirect product of $\C^*$ and
$\C_+$ acts on $V$ is equivalent to the condition that $\p$ is a
homogeneous derivation (cf.\ Lemma \ref{sdp}). Thus we are led to
pairs
$$
(A,\p),\qquad e=\deg \p,
$$
where $\p$ is a homogeneous locally nilpotent derivation on $A$ of
a certain degree $e$. Our classification of such pairs is as follows.

\smallskip

{\em Elliptic case}: In this case $A_0\cong \C$, and $A$ is
positively graded so that the associated $\C^*$-surface $V=\Spec
A$ has a unique fixed point given by the maximal ideal
$A_+:=\bigoplus_{i>0}A_i$. If $V$ also admits a nontrivial
$\C_+$-action then by \cite[Lemmas 2.6 and 2.16]{FlZa2},
$V\cong\A^2_\C/\Z_d$ is a quotient of $\A^2_\C$ by a small cyclic
subgroup of $\GL_2(\C)$. More precisely, we show in Theorem
\ref{linact} that {\em $A\cong \C[X,Y]^{\Z_d}$, where the cyclic
group $\Z_d:=\Z/d\Z=\langle \zeta\rangle$ generated by a
primitive $d$-th root of unity $\zeta$ acts on $\C[X,Y]$ via
$\zeta . X=\zeta X$ and $\zeta . Y = \zeta^e Y$ with $e\ge
0,\,\,\,\gcd(e,d)=1$, and $\p=X^e\p/\p Y$.}  In particular,
$V\cong V_{d,e}$ is an affine toric surface (see Example
\ref{tosu}).

\smallskip

{\em Parabolic case}: Here again $A$ is positively graded,
but
$A_0\ne \C$. Thus $C=\Spec A_0$ is a smooth affine curve,
and $V$
is fibered over $C$ with general fiber $\A^1_\C$. Using the
DPD construction  it follows that $A=A_0[D]$
for some
$\Q$-divisor $D$ on $C$ (see \cite[Theorem 3.2]{FlZa1}).
More
precisely, if $K_0$ denotes the field of fractions
$\Frac(A_0)$
then $A=A_0[D]\subseteq K_0[u]$ is the graded subring with
$$
A_n=\{fu^n\in K_0\cdot u^n\, | \,\div\,f+nD\ge 0\}\,.
$$
If such a surface admits also a $\C_+$-action given by a
homogeneous locally nilpotent derivation $\p$ then either $\C_+$
acts {\em vertically} (that is fiberwise), so that the orbits are
contained in the fibers of the projection $V\to C$, or the orbits
map onto the base curve $C$ ({\em horizontal case}). In both
cases we  classify all possible actions (see Theorems
\ref{strthm} and \ref{TN}). For instance, in the horizontal case
$V\cong V_{d,e}\cong \A_\C^2/\Z_d$ is again an affine toric
surface and the derivation $\p$ is as described in the elliptic
case. These are the only normal affine surfaces with an elliptic
or parabolic $\C^*$-action and with a trivial Makar-Limanov
invariant that is, admitting two non-trivial $\C_+$-actions with
different orbit maps (see Definition \ref{mli} and Theorem
\ref{zigzag}).

\smallskip

{\em Hyperbolic case}: In this case $A_i\ne 0$ for all $i\in \Z$,
and the surface $V=\Spec A$ is fibered over the base curve
$C=\Spec A_0$ with general fiber $\C^*$. By \cite[Theorem
4.3]{FlZa1} $A=A_0[D_+,D_-]\subseteq \Frac(A_0)[u,u^{-1}]$ with a
pair of $\Q$-divisors $D_\pm$ on $C$ satisfying $D_++D_-\le 0$.
This means that $A_{\ge 0}=A_0[D_+]\subseteq K_0[u]$ and $A_{\le
0}=A_0[D_-]\subseteq K_0[v]$ are as above with $v=u^{-1}$.
Furthermore, the pair $(D_+,D_-)$ is determined uniquely up to an
arbitrary shift $(D_+,D_-)\rightsquigarrow (D_++\div\,
\varphi,D_--\div\, \varphi)$ with $\varphi\in \Frac A_0$. In
Corollary \ref{crucol} we show that $A$ {\em admits a homogeneous
locally nilpotent derivation $\p$ of positive degree $e$ if
and only if $C\cong \A^1_\C$ i.e., $A_0\cong \C[t]$, and
$A\cong A_0[D_+,D_-]$, where $D_+=-\frac{e'}{d}[p]$ is
supported at one point, $\,0\le e'<d$ and $ee'\equiv 1 \mod
d$. Moreover, $\p$ is uniquely determined up to a constant by
its degree.} Alternatively, such surfaces can be described as
cyclic quotients of the normalizations of hypersurfaces
$\{u^dv-p(t)=0\}$ in
$\A_\C^3$, where $p\in\C[t]$ (see \cite[Proposition 4.14]{FlZa1}
and Corollary \ref{maincor} below).

$\C_+$-actions on a normal affine surface $V$ are related to
affine rulings $V\to \Gamma$   (that is, rulings into affine
lines) with $\Gamma$ being a smooth affine curve (see Lemma
\ref{afru}). If $V=\Spec A$ with $A=A_0[D_+,D_-]$ as above, where
$A_0=\C[t]$ and $D_++D_-\ne 0$, then {\it there exists an affine
ruling $V\to \A_\C^1$ if and only if the fractional part
$\{D_\pm\}$ of at least one of the $\Q$-divisors $D_\pm$ is
supported at one point or is zero. Such an affine ruling is unique
unless both $\{D_+\}$ and $\{D_-\}$  are either zero or supported
at points $\{p_\pm\}$, and if and only if, for a homogeneous
element $v\in A\backslash \C$, $\ker \p\supseteq \C[v]$ for every
locally nilpotent derivation $\p\in\Der A$} (Corollary
\ref{crucol} and Theorem \ref{long}). Otherwise $V$ allows
continuous families of affine rulings, of $\C^*$-actions and of
$\C_+$-actions with generically different orbit maps (Corollary
\ref{cf}). The same is also true in the elliptic and the
parabolic cases.

\smallskip

In the first two sections we summarize some facts on
$\C_+$-actions and on algebraic group actions on normal affine
surfaces. Section 3 contains the principal classification
results. In section 4 we classify all $\C^*$-surfaces which have
a trivial Makar-Limanov invariant (Corollary \ref{epml} and
Theorem \ref{long}). Finally, in Section 5 we discuss concrete
examples and compare different approaches.

Throughout the paper we use the notation
$\GL_2=\GL(2,\C),\,\,\,\SL_2=\SL(2,\C)$, etc.

\section{$\C_+$-actions
and locally
nilpotent derivations}
We frequently use the following well
known facts.

\bprop\label{Ren}{\em (see e.g., \cite{Ren, ML1, Za})} Let
$V=\Spec A$ be an affine algebraic $\C$-scheme. Then the
following hold:

\begin{enumerate}
\item[(a)] If $\C_+$ acts on $V$
then the associated derivation $\partial$ on $A$ is locally
nilpotent, i.e.\ for every $f\in A$ we can find $n\in \N$ such
that $\partial^n(f)=0$. Conversely, given a locally nilpotent
$\C$-linear derivation $\partial:A\to A$ the map
$\varphi:\C_+\times A\to A$ with $\varphi(t,f):=e^{t\partial}f $
defines an action of $\C_+$ on $V$.

\item[(b)] Assume that $A$ is a domain and
let $\p\in {\Der}_{\C} A$ be a locally nilpotent derivation of
$A$. Then the subalgebra $\ker \p=A^{\C_+}\subseteq A$ is
algebraically and factorially closed (or inert)\footnote{The
latter means that $ab\in\ker \p\Rightarrow a,\,b\in\ker \p$.} in
$A$, and the field extension $\Frac (\ker \p) \subseteq \Frac A$
has transcendence degree $1$. Moreover, for any
$u\in\Frac A$ with
$u\p(A)\subseteq A$, the derivation $u\p\in \Der_{\C} A$ is
locally nilpotent if and only if $u\in\Frac(\ker\p)$.

\item[(c)] If $\C_+$ acts non-trivially
on an irreducible reduced affine curve $C$ then $C\cong\A^1_\C$.
\end{enumerate}
\eprop

\bcor\label{norm} For an algebraic $\C$-scheme $A$ and a
locally nilpotent derivation $\p$ on $A$, the
following hold.
\begin{itemize}
\item[(a)] The
algebra of invariants
$\ker\p=A^{\C_+}$ is integrally closed in $A$. Consequently, if
$A$ is normal and the ring of invariants $A^{\C_+}$ is finitely
generated then the orbit space $\Spec\,A^{\C_+}$ of the
associate $\C_+$-action on $V$ is also normal.

\item[(b)] For an element $v\in A$, the principal ideal
$(v)=vA$ is $\p$-invariant if and only if $v\in \ker \p$.

\item[(c)] If $\dim A\ge 2$ then the automorphism group $\Aut\,A$ is
of infinite dimension. \end{itemize} \ecor

\bproof (a) immediately follows from Proposition \ref{Ren}(b). To
show (b) we fix $n\ge 1$ such that $u:=\p^{n-1}(v)\neq 0$ and $\p
u=0$. If the ideal $(v)$ is $\p$-invariant then $u\in \ker \p \cap
(v)$ can be written as $u=fv$ with $f\in A$. As $\ker \p$ is inert
(see Proposition \ref{Ren}(b)) and $u=fv\in \ker \p$ we have $v\in
\ker \p$, as required. The proof of the converse is trivial. As
$e^{a \p}\in \Aut\,A$ $\forall a \in \ker\p$ and $\dim \ker
\p \ge 1$, (c) also follows from Proposition \ref{Ren}(b).
\eproof

\bsit\label{cld} Let us recall some well known facts on the
surface geometry in presence of a $\C_+$-action; see e.g.,
\cite{BaML2, Miy2, MaMiy}. For a normal affine surface $V$ we
denote $V_{\rm reg}=V\backslash \Sing V$. A {\it cylinder} in $V$
is a Zariski open subset $U\cong \Gamma_0\times\A^1_\C$, where
$\Gamma_0$ is a smooth curve. An {\it affine ruling} on $V$ is a
morphism $V\to \Gamma$ onto a smooth curve $\Gamma$ with general
fibers isomorphic to $\A^1_\C$. Two affine rulings coincide if
they have the same fibers. \esit

\blem\label{afcy} (\cite[Ch. 3, Lemma 1.3.1,
Theorem 1.3.2 and Lemma
1.4.4(1)]{Miy2}) For a normal affine surface
$V$ the following conditions are equivalent:
\begin{enumerate}
\item[(i)] $V$ is affine ruled.
\item[(ii)] $V$ contains a cylinder.
\item[(iii)] There exists an affine Zariski open subset $W\subseteq
V_{\rm reg}$ with ${\bar k}(W)=-\infty$. \footnote{As usual,
${\bar k}$ stands for the logarithmic Kodaira dimension.}
\end{enumerate}

\noindent Moreover, under these conditions $V$ has at most cyclic
quotient singularities. \elem

\brem\label{equi}
If $V$ is smooth then any degenerate fiber of
an affine ruling on $V$ consists of disjoint components
isomorphic to $\A^1_\C$ (see \cite{Be1, Fi}). If $V$ is only
normal then any such component has a normalization isomorphic to
$\A^1_\C$, contains at most one singular point of $V$ and is
smooth off this point (\cite[Ch. 3, Lemmas 1.4.2 and
1.4.4]{Miy2}).
\erem

Suppose that a normal surface $V=\Spec \,A$  admits a non-trivial
$\C_+$-action. The orbit morphism $\pi_+:V\to \Gamma:=V\quot\C_+$
then yields an affine ruling on $V$ over a smooth affine curve
$\Gamma\cong \Spec\,A^{\C_+}$. Therefore \cite[Remark 1]{BaML2},
an affine ruling on $V$ over a projective base cannot be produced
in this way. For instance, the latter concerns the projection
pr$_1 : (\pP^1\times\pP^1)\backslash \Delta \to \pP^1$, where
$\Delta \subseteq \pP^1\times\pP^1$ is the diagonal. The
following simple lemma clarifies the situation (cf. \cite[Prop.
2]{BaML2}).

\blem\label{afru}
For a normal affine surface $V$ the following
are equivalent:
\begin{enumerate}
\item[(i$'$)]
$V$ admits an affine ruling $V\to \Gamma$ over an affine base
$\Gamma$.
\item[(ii$'$)] $V$ contains a cylinder
$U\cong \Gamma_0\times\A^1_\C$ which is
a principal Zariski open subset.
\item[(iii$'$)] There exists a non-trivial regular
$\C_+$-action on $V$.
\end{enumerate}\elem

\bproof The implication (iii$'$) $\Rightarrow$ (i$'$) has been
noted above. The proof of (i$'$) $\Rightarrow$ (ii$'$) follows
that of (i) $\Rightarrow$ (ii) in Lemma \ref{afcy}; it suffices
to note that, because $\Gamma_0\subseteq \Gamma$ can be taken
principal, so does the cylinder $U\subseteq \pi^{-1}(\Gamma_0)$.

To show (ii$'$) $\Rightarrow$ (iii$'$) we let $U\cong
\Gamma_0\times\A^1_\C$ be a principal cylinder in $V=\Spec\,A$
given via $A[1/f_0]\cong B[t]$ with $f_0\in A$, where
$\Gamma_0=\Spec B$. We consider the derivation $\p=\p/\p t\in
\Der\,B[t]$. Given a system of generators $g_1,\ldots,g_n$ of the
algebra $A$ we can write $\p g_i={h_i\over f_0^{k_i}}$, where
$h_i\in A$ and $k_i\ge 0\quad (i=1,\ldots,n)$. Since $f_0^e\p
(g_i)\in A\quad\forall i$, where $e:=\max_{1\le i\le n} \,k_i$, we
have $\p_e:=f_0^e\p \in \Der\,A$. Moreover, $\p_e f_0=0$ as
$f_0\in B[t]$ is a unit. Hence $\p_e$ is locally nilpotent and
so defines a non-trivial $\C_+$-action on $V=\Spec \,A$, as
required. \eproof

\bsit\label{sti} If a ramified covering of normal varieties $Y\to
X$ is unramified in codimension $1$ then any $\C_+$-action on $X$
lifts to $Y$ \cite[proofs of Lemmas 2.15 and 2.16]{FlZa2}. In the
following lemma we show that, under certain circumstances, it
still lifts to a cyclic covering ramified in codimension 1,
provided the latter is defined by an invariant.\esit

\blem\label{cyex}
Let $A$ be a normal domain of finite type over
$\C$ and let $\p \in \Der A$ be a non-zero locally nilpotent
derivation. For a non-zero element $v\in \ker \p$  and for $n\in
\N$ denote $A'$ the normalization of the cyclic ring
extension
$A[u']\supseteq A$ with $(u')^n=v$. Then the following hold:

\begin{enumerate}
\item[(a)] $A'$ is a normal affine $\C$-algebra of finite type, and
the elements of $A$ are not zero divisors on $A'$.

\item[(b)] $\p $ extends uniquely to a locally nilpotent
derivation $\p'\in\Der A'$ with $\p' (u')=0$.

\item[(c)] If, moreover, $A$ is a graded domain and $v$ and $\p$
are homogeneous with $\deg v =n$ then $A'$ is graded as well, and
$u'$ and $\p'$ are homogeneous with $\deg u'=1$ and $\deg\p'=\deg\p$.

\item[(d)] Furthermore, if the polynomial $x^n-v\in A[x]$
is irreducible over
$A$ then the cyclic group $\Z_n=\langle \zeta\rangle$, where
$\zeta$ is a primitive $n$-th root of unity, acts on $A'$ with
$\zeta | A={\rm id},\quad\zeta . u'=\zeta u'$, and $A=(A')^{\Z_n}$
is the ring of invariants of this action.

\end{enumerate}\elem

\bproof The proofs of (a), (c) and (d) are easy and we omit
them.
To show (b) note that the derivation
$\p'\in\Der A[u']$ defined by $\p' |A=\p$ and $\p'(u')=0$ is
locally nilpotent. By \cite{Sei} its extension to $\Frac (A[u'])$
stabilizes the integral closure $A'$ of $A[u]$. By \cite{Vas}
(see also \cite[Lemma 2.15(a)]{FlZa2}) this extension $\p'$ of
$\p$ to $A'$ is again locally nilpotent, as stated. \eproof

\section{Algebraic group actions on affine
surfaces}

\subsection{$\C_+$-actions on graded rings} We let $V=\Spec A$
be an affine variety over $\C$
with an effective $\C^*$-action,
which corresponds to a grading $A=\bigoplus_{i\in\Z} A_i$.

\blem\label{rl} {\em \cite{Ren}} If $\partial$ is a locally
nilpotent derivation on $A$ and $\partial=\sum_{i=k}^l\partial_i$
is the decomposition of $\partial$ into graded components then
$\delta_k$ and $\delta_l$ are again locally nilpotent. \elem

     Homogeneous locally nilpotent
derivations on $A=\bigoplus_{i\in\Z} A_i$ correspond to actions
of certain semidirect products of $\C^*$ and $\C_+$ on $A$.
Indeed, we have the following lemma (cf. \cite{Po},
\cite[(2.5)]{Be2}).

\blem \label{sdp} \begin{enumerate}\item[(a)] Let
$\partial:A\to A$ be a homogeneous locally nilpotent derivation
of degree $e$ and consider the action of $\C^*$ on $\C_+$
given by $\tau_e(t,\alpha):=t.\alpha=t^e\alpha$, where $t\in\C^*$,
$\alpha\in\C_+$. Then the semidirect product
$$
G_e:=\C^* \ltimes_{\tau_e}\C_+
$$
(with $\C_+$ as a normal subgroup) acts on $A$, and hence on
$V$, via
$$
(s,\alpha).f:= s. e^{\alpha\partial}(f),\quad \text{where }
(s,\alpha)\in G_e \text{ and } f\in A.
$$
This action restricts to the given actions on the subgroups
$\C_+$ and $\C^*$ of $G_e$.
\item[(b)] Conversely, if
there is an action of $G_e$ on $V=\Spec A$ restricting to the
given action of $\C^*\subseteq G_e$ on $A$, then $\C_+\subseteq
G_e$ acts on $V$ and the associated derivation $\partial$ on
$A$ is homogeneous of degree $e$. \end{enumerate}
\elem

\bproof (a) The multiplication on $G_e$ is given by
$$(s,\alpha)(t,\beta)=(st,t^{-e}\alpha+\beta)\qquad
{\rm with}\quad s,t\in \C^*,\,\,\,\alpha,\beta\in\C_+\,.$$ Since
$\partial$ is homogeneous of degree $e$ it follows that
$s.(\partial(f))=s^e\partial(s.f)$, and so
$$s.e^{\alpha\partial} (f)=\sum_{\nu=0}^{\infty}
s.{\alpha^{\nu}\partial^{\nu}(f)\over\nu!}
=\sum_{\nu=0}^{\infty}{\alpha^{\nu}s^{e\nu}
\partial^{\nu}(s.f)\over\nu!}
=e^{s^e\alpha\partial}(s.f)\,,$$ hence
$$\left((s,\alpha)(t,\beta)\right).f=(st).e^{(t^{-e}
\alpha+\beta)\partial}(f)= s.(t.e^{t^{-e}\alpha\partial}
e^{\beta\partial}(f))= s.e^{\alpha\partial} t.e^{\beta\partial}(f)
=(s,\alpha).((t,\beta).f)\,.$$ This shows that $G_e$ acts indeed
on $A$ and hence on $V$.

(b) Conversely, suppose that $G_e$ acts on $A$ restricting to the
given action of $\C^*$ on $A$. Then for $\C_+\ni a= (1,a)\in G_e$
we have $a.f=e^{a\partial}(f)$, and so
$$s.e^{a\partial}(f)=(s,a).f=(1,s^ea)((s,0).f)=
e^{s^ea\partial}(s.f)\,.$$ Differentiating this equation with
respect to $a$ and taking $a=0$ one gets
$$s.\partial(f)=s^e\partial(s.f)\,.$$ It follows that
$\partial$ is homogeneous of degree $e$. \eproof

\brems\label{newder} 1. For any non-zero homogeneous element $u\in
\ker\p$ of degree $n$, the derivation $\p':=u^m\p\in\Der\,A$ is
again locally nilpotent (see Proposition \ref{Ren}(b)) of degree
$e+mn$. Thus for every $m\ge 0$ the group $G_{e+mn}$ also acts on
$A$ restricting to the given $\C^*$-action on $A$. The inversion
$\lambda\longmapsto \lambda^{-1}$ provides an isomorphism
$G_e\cong G_{-e}$, and so $G_{e'}$ acts on $V$ for
any $e'\equiv \pm e\mod n$.

\no 2. For instance, a Borel subgroup $B\subseteq \SL_2$ is
isomorphic to $G_2$ and  acts effectively on $V=\A^2_C$ with an
open orbit. Similarly, the Borel subgroup $B':=B/\Z_2$ in
$\PGL_2=\SL_2/\Z_2$, where $\Z_2=\{\pm I_2\}$ is the center of
$\SL_2$ (and of $B$), is isomorphic to $G_1$ and acts effectively
on the Veronese cone $V_{2,1}:=\A_\C^2/\Z_2\cong \Spec
\C[t,u,v]/(uv-t^2)\subseteq \A_\C^3$ with an open orbit (cf.
Example \ref{pop}).

  \no 3. For $e>0$, $G_e$ is a metabelian solvable Lie group
  with a cyclic center
  $Z(G_e)=\Z_e\ltimes\{0\}\subseteq \C^*\ltimes \C_+$, and so
  is an \'etale covering group of $G_1$ via
  $G_e\stackrel{e:1}{\to} G_1\cong G_e/Z(G_e)$.
  The Lie algebra
  $\fg=$Lie$\, G_e$ is isomorphic to $\A_\C^2$ with Lie bracket
  $[\vec v_1,\vec v_2]=(0,\vec v_1\wedge \vec v_2)$. \erems

Actually, an effective $G_e$-action on $A$ with $e\neq 0$ permits to
produce a continuous family of gradings on $A$.

\bprop\label{roun} Let $A=\bigoplus_{i\in\Z} A_i$ be a graded
$\C$-algebra of finite type and $\p\in\Der\,A$ be a homogeneous
locally nilpotent derivation on $A$ of degree $e\neq 0$. If the
orbit closures of the associated $\C^*$- and $\C_+$-actions on
$V:=\Spec A$ are generically different then $A$ admits a
continuous family of generically distinct gradings.\eprop

\bproof For $\alpha\in\C_+$, $\alpha\neq 0$, we consider a new
$G_e$-action on $V_\alpha:=V$ induced by the isomorphism $\alpha
: V \to V_\alpha$ that is, conjugated with the original
$G_e$-action on $V$ by means of $\alpha$. More precisely, we have
a commutative diagram
\begin{diagram}
G_e & \rTo^{\xi_\alpha} & G_e\\
\dTo<{-.x} && \dTo>{-.x}\\
V& \rTo^{\alpha.-} & V_\alpha
\end{diagram}
where the vertical arrow on the right is the new $G_e$-action on
$V_\alpha$ and $$\xi_\alpha (g) = \alpha g \alpha^{-1}=(t, \,\beta
+t^{-e}\alpha-\alpha)\qquad\mbox
{for}\qquad g=(t, \,\beta)\in G_e\,.$$
The
$\C^*$-orbit of
$(1,\,\beta)\in G_e=\C^* \ltimes_{\tau_e}\C_+$ is equal to $\C^*
\times \{\beta\}$ and is mapped under $\xi_\alpha$ onto the set
$$\left\{(t, \,\beta+t^{-e}\alpha-\alpha) \mid t \in
\C^*\right\}\,,$$ which is not an orbit of the $\C^*$-action on
$G_e$. Since by our assumption for a general $x\in V$ the orbit
$G_e . x$ has dimension 2, the generic $\C^*$-orbit in $V$ is not
mapped onto a $\C^*$-orbit of $V_\alpha$. \eproof

In the surface case we have the following elementary lemma.

\blem\label{quo} For a $G_e$-action on an affine surface
$V=\Spec\,A$ the following conditions are equivalent.

\begin{enumerate}
\item[(i)] It has an open orbit.

\item[(ii)] $A^{\C_+}\neq A^{\C^*}$ ($\Leftrightarrow \ker\p\neq A_0$)
\footnote{I.e., the $\C_+$-action is horizontal w.r.t. the given
$\C^*$-action.}.

\item[(iii)] $\ker\p=\C[v]$ or $\ker\p=\C[v,v^{-1}]$, where $v\in
A_d$ with $d\neq 0$.

\end{enumerate} Under these equivalent conditions the surface
$V$ is rational, and the affine ruling $v: V\to
\Gamma:=\Spec\,A^{\C_+}$ has at most one degenerate fiber $v=0$
consisting of $\C^*$-orbit closures \footnote{Cf. Remarks
\ref{equi}, \ref{adlnd}.iii and Lemma \ref{trmap} below.}.
\elem

\bproof Since $\p\in\Der\,A$ is homogeneous, its ring of invariants
$\ker\p=A^{\C_+}$ is a graded subring of $A$.  Thus the normal
(hence smooth) affine curve $\Gamma=\Spec\,A^{\C_+}$ also carries
a $\C^*$-action, and the quotient morphism $V\to \Gamma=V\quot\C_+$
(which provides an affine ruling on $V=\Spec\,A$) is
$\C^*$-equivariant. In case $A^{\C_+}\neq A^{\C^*}$
(that is, $\ker\,\p\neq A_0$) the induced $\C^*$-action on
$\Gamma$ is non-trivial, hence $\Gamma\cong\A^1_\C$ or $\C^*$. In
this case $\ker\,\p=\C[v]$ and $\C[v,v^{-1}]$, respectively,
where $v\in A_d \cap \ker\,\p$ is homogeneous and $d\neq 0$.

The rationality of $V$ follows from L\"uroth's Theorem. The rest
of the proof is easy and can be omitted.\eproof

\subsection{Actions with an open orbit}
The next simple observations will be used in the proofs below (cf.
Remark 1 in \cite[II.4.3.B]{Kr}).

\blem\label{obser} \begin{enumerate} \item[(a)] If a connected Lie
group $L$ and a finite group $G$ act on an algebraic $\C$-scheme
$V=\Spec A$ then the action of $L$ descends to $V/G$ if and only
if the actions of $G$ and $L$ on $V$ commute.

\item[(b)] Conversely, suppose that a connected and simply
connected Lie group $L$ acts on the quotient $V/G$ of $V$ by a
free action of a finite group $G$. Then the action of $L$ lifts
to $V$ commuting with the action of $G$. \end{enumerate}\elem

\bproof (a) Suppose first that the action of $L$  on $V$ descends
to $V/G$. We may also assume that $G$ acts faithfully on $V$. It
follows that $L$ preserves the $G$-orbits, and so, if $w=g.z$ for
some $z,w\in V$ and some $g\in G$ then for any $\lambda\in L$
there is an element $g'=g'(\lambda)\in G$ such that $\lambda .
w=g' . \lambda . z$. This implies the equality $\lambda
g.z=g'\lambda .z$. Since $g'(\lambda)$ is a continuous function
on the connected Lie group $L$ with values in $G$ it must be
constant, i.e., $g'=g$, and so $g\lambda=\lambda g $ for all
$g\in G$ and $\lambda\in L$. Thus the actions of $L$ and of $G$
commute, as stated  in (a). The  proof of the remaining
assertions is easy and will be omitted.
\eproof

\blem\label{gtl} \begin{enumerate}
\item[(a)] If a complex
unipotent Lie group $U$ acts on an affine variety $V$
with an open orbit then $V\cong \A_\C^{\dim V}$.

\item[(b)] If a complex reductive Lie group $G$
acts effectively on a connected algebraic variety $V$ with a
fixed point
$p\in V$ then the induced representation $\tau_p : G \to
\GL (T_p V)$ on the Zariski tangent space of $V$ at $p$
is faithful.

\item[(c)] Any affine toric surface $V$
non-isomorphic to $\C^*\times\C^*$ admits  a $G_l$-action with an
open orbit for every $l\in\Z$.
\end{enumerate}
\elem

\bproof (a) Since any orbit of $U$ is closed in $V$
\cite[Exercise 8 in Section 17]{Hu1}, \cite[III.2.5.3]{Kr}, the open
$U$-orbit is the whole $V$. Thus $V\cong U/H\cong\A_\C^{\dim V}$,
where
$H\subseteq U$ is a closed subgroup (see \cite[Corollary of Theorem
2]{Po}). This shows (a).

(b) is well known and follows for instance from Luna's \'etale
slice theorem or from the identity theorem
\cite[Sect.\ 2.1]{Ak}. Alternatively, this can be seen
by the following elementary argument: for $n\gg0$ the induced
action of
$G$ on
$A_n:=\cO_{V, p}/\fm^{n+1}$ is easily seen to be faithful,
i.e.\ the map $\rho_n: G\to \Aut(A_n)$ is injective, where
$\Aut(A_n)$ denotes the Lie group of $\C$-algebra automorphisms of
$A_n$. The subgroup $N_n$ of $\Aut(A_n)$ consisting of automorphisms
$f$ with $f\equiv \id \mod \bar\fm^2$ is a normal unipotent
subgroup, so
$\rho_n^{-1}(N_n)$ is also normal and unipotent and thus
trivial. It follows that already the map $G\to
\Aut(A_1)\cong \Aut(A_n)/N_n$ is injective, which implies that
$G$ acts effectively on $T_pV$.

(c) As $G_l\cong \C^*\times\A_\C^1$ this
is evident in case that $V\cong \C^*\times\A_\C^1$.
Otherwise (c) is shown in Example
\ref{tosu} below.
\eproof

\bexa\label{tosu} {\it Affine toric surfaces.} Given two natural
numbers $d,\, e'$ with $0\le e'<d,\,\,\gcd (e',d)=1$, we consider
the affine toric surface $V_{d,e'}=\Spec A_{d,e'}$, where
\be\label{tor1} A_{d,e'}=\C[X,Y]^{\Z_d}\cong \bigoplus_{b\ge
0,\,\, ad-be'\ge 0} \C \cdot x^{a}y^{b}\subseteq \C[x, y] \ee is
the semigroup algebra of the cone $\sigma^{\vee}=C\,({\vec e}_1,
e' {\vec e}_1+d {\vec e}_2)$ in $\R^2$, and where
$\Z_d=\langle\zeta\rangle$ acts on $\C[X,Y]$ via
\be\label{zd2}\zeta .X=\zeta X,\qquad\zeta .Y=\zeta^{e'} Y\,\ee
(cf. \cite[Example 2.3]{FlZa1}). This $\Z_d$-action commutes with
any $\C^*$-action on $\C[X,Y]$ of the form
$$
\lambda .X=\lambda^{d_X} X,\qquad \lambda . Y=\lambda^{d_Y} Y\,,
$$ where
$(d_X,d_Y)\in\Z^2$. It also commutes with the locally nilpotent
derivations
\be \label{torder}
\p_{X, e''}=X^{e''}{\p\over\p Y}\qquad\mbox{and}\qquad \p_{Y,e}
=Y^{e}{\p\over\p X}\quad\in \Der\ \C[X,Y]\,, \ee
where $ e''$, $e\ge 0$ are such that $e''\equiv e'\mod d$ and
$e\cdot e'\equiv 1\mod d$ if
$e'\neq 0$, $e=0$ if $e'=0$. Therefore by Lemma \ref{obser} the
$\C_+$-actions on $\C[X,Y]$ induced by $\p_{X, e''}$ and
$\p_{Y,e}$ stabilize the ring of $\Z_d$-invariants
$A_{d,e'}=\C[X,Y]^{\Z_d}$, hence descend from $\A_\C^2=\Spec
\C[X,Y]$ to the quotient surface $V_{d,e'}=\Spec
A_{d,e'}=\A_\C^2/\Z_d$. Note that any affine toric surface
non-isomorphic to $\C^*\times\C^*$ or $\A^1_\C\times\C^*$, is
isomorphic to $V_{d,e'}$ for some $d,e'$ as above. Consequently,
any such surface admits two $\C_+$-actions with different general
orbits (cf. Corollary \ref{epml} below).

Letting above e.g., $d_X=0,\,\, d_Y=-l$ we obtain that $\deg
\p_{X, e'}=l$, and so by Lemma \ref{sdp}(b) the group $G_l$ acts
effectively on the ring $A=A_{d,e'}$.\eexa

\blem\label{EE} Let $G$ be a connected complex algebraic Lie
group acting effectively on a normal affine surface $V=\Spec A$.
\begin{enumerate}
\item[(a)] If $G$ is unipotent and
  $V\not\cong\A_\C^2$, then $G$ is commutative and the orbits of
$G$ are 1-dimensional.
  \item[(b)] If $G$ is solvable and acts on $V$ with an open
orbit $O$, then $O$ is isomorphic to one of the surfaces
$\C^*\times\C^*$, $\C^*\times\A_\C^1$ or $\A_\C^2$. Moreover, if
$O$ is big that is, $V\backslash O$ is finite, then $O=V$.
\item[(c)] $G$ is solvable if and only if it does not contain a
subgroup isomorphic to $\SL_2$ or to $\PSL_2$.
\end{enumerate}
\elem

\bproof (a) The orbits of $G$ are closed in $V$ and generically
one-dimensional, since otherwise $V\cong \A_\C^2$ by Lemma
\ref{gtl}(a). We let $\pi:V\to \Gamma:= \Spec A^G$ be the
quotient map. The Lie algebra $\fg=$Lie$\,G$ consists of vector
fields tangent along the fibers of $\pi$. Any such vector field
$\p\in\fg$ is an infinitesimal generator of a one-parameter
subgroup of $G$ isomorphic to $\C_+$ and so is a locally
nilpotent derivation on $A$. Being proportional, every two such
nonzero derivations $\p_1,\p_2$ are equivalent i.e.,
$b_1\p_1=b_2\p_2$ for some $b_1,b_2 \in A^G$. Thus $\p_1=b\p_2$
with $b:=b_2/b_1\in \Frac A^G$ and so
$0=[\p_1,b\p_2]=b[\p_1,\p_2]$. This shows that $\p_1$ and $\p_2$
commute, proving (a).

(b) We may suppose that $V\not\cong\A_\C^2$. In the decomposition
$G=\T \ltimes N$ \cite[Theorem 19.3(b)]{Hu1}, where $\T$ is a
maximal torus and  $N$ is the unipotent radical of $G$, we have
$N\cong\C_+^r$ by (a). If $r=0$ then clearly $V\cong
\C^*\times\C^*$. In case $r>0$ let $\p_0\in \,$Lie$\,N$ be a
common eigenvector of the adjoint representation of $\T$ on
Lie$\,N$ and denote $N_0\subseteq N$ the corresponding
one-parameter subgroup. By (a) the orbits of $G$ and of $G_0:=\T
\ltimes N_0$ are the same. Thus we may suppose that $N=N_0$ has
dimension 1. As $G$ acts effectively on $V$ with an open orbit
the torus $\T$ must be of dimension 1 or 2, so
$G\cong \C^*\ltimes \C_+$ or $G\cong {\C^*}^2\ltimes \C_+$. In
the first case the open orbit $O$ of $G$ is isomorphic to $G$.
In case $G\cong {\C^*}^2\ltimes \C_+$ the stabilizer
$H=$Stab$_x\subseteq G$ of a point $x\in O$ has dimension 1 and
so $H=N$ or $H\cong\C^*$. If $H=N$ then $O\cong
G/H\cong{\C^*}^2$.  If $H\cong\C^*$ then we may suppose that
$H\subseteq \T$.  Indeed, any subtorus in $G$ is contained in a
maximal torus, which is unique up to a conjugation. But then
$O\cong G/H\cong \C^*\times\A_\C^1$.

In all cases the open orbit $O$ is affine, hence $V\backslash O$
is either empty or a divisor. Thus, if $O$ is big then $O=V$,
proving (b).

(c) is well known and follows from the structure theory of
algebraic groups, see
\cite{Bou, Hu1}.

%
\eproof

To describe  all normal affine surfaces $V$ admitting an
action of an algebraic group $G$ with an open (not necessarily
big) orbit, we follow a suggestion in \cite[The concluding
remark]{Po}. In the particular case of smooth rational surfaces
it was confirmed in \cite[Proposition 2.5]{Be2}.

\bprop\label{agoo}
Let $V=\Spec\,A$ be a normal affine
surface non-isomorphic to
$\C^*\times\C^*$. If an
algebraic group $G$ acts on $V$ with an open orbit then, for some
$e\in\Z$, the group $G_e=\C^*
\ltimes_{\tau_e}\C_+$ also acts on $V$ with an open orbit.
\eprop

\bproof
If $V$ is a toric surface then by Lemma
\ref{gtl}(c) it admits a $G_e$-action with an open orbit. So we
may suppose in the sequel that $V$ is not toric.

In case $G\cong\SL_2$  we let $B_\pm$ be the Borel subgroups of
upper/lower triangular matrices. Their intersection is the torus
$\T\cong\C^*$ of diagonal matrices. If both $B_\pm$ act with
1-dimensional orbits on $V$ then their orbits would be equal to
the orbit closures of the torus action. Hence also $G$ would act
with 1-dimensional orbits contradicting our assumption. Thus at
least one of the groups $B_\pm$ has an open orbit in $V$. Since
$B_\pm\cong G_2$  the result  follows in this case.

  Clearly, the case $G\cong \PGL_2\cong\SL_2/\{\pm I\}$
reduces to the previous one.

  For the remaining cases we may suppose that $G$ acts
effectively on
$V$, is connected and does not contain a subgroup isomorphic to
$\SL_2$ or $\PGL_2$. By Lemma \ref{EE}(a),(c) $G$ is solvable and
not unipotent. As $V$ is not toric, the maximal torus $\T$ of $G$
has dimension 1. As in the proof of Lemma \ref{EE}(b) we can
  restrict the action of $G$ to a subgroup $H=\T\ltimes \C_+$
of $G$ which still has an open orbit. As $H\cong G_e$ for some
$e$, the result follows. \eproof

\section{Classification of  affine surfaces with a $\C^*$-
and $\C_+$-action}

In this section we study normal affine surfaces $V=\Spec A$
endowed with an effective $\C^*$- and a $\C_+$-action. The
$\C^*$-action provides a grading $A=\bigoplus_{i\in\Z}A_i$ and
the $\C_+$-action a locally nilpotent derivation $\p$ of $A$. Due
to Lemma \ref{rl} we can  find a {\em homogeneous} locally
nilpotent derivation on $A$. Thus  in the sequel we consider
pairs $(A,\p)$,  where $A$ is the graded coordinate ring of
$V=\Spec A$ as above and $\p\in\Der A$ is a nonzero homogeneous
locally nilpotent derivation.

\bdefi\label{ellpa}We call such a pair $(A,\p)$ {\em elliptic}
if the $\C^*$-action on $V$ is elliptic i.e., if $A$ is
positively graded with\
$\dim A_0=0$, {\it
parabolic} if $A$ is parabolic i.e., positively
graded with\
$\dim A_0=1$, and {\it
hyperbolic} if
$A$ is hyperbolic,  i.e.\ $A_\pm\neq 0$.

Two such pairs $(A,\p)$ and $(A',\p')$ are called {\it isomorphic}
if there is an isomorphism of graded $\C$-algebras $\varphi:A\to
A'$ with $\varphi \p=\p' \varphi$.

For hyperbolic pairs we will suppose in the sequel that
$e:=\deg\p\ge 0$ (indeed, otherwise we can reverse the grading of
$A$).\edefi

We can reformulate \ref{sdp} in this setup as follows.

\bprop\label{groupder}
Let $e\in \Z$ be fixed.
There is a 1-1 correspondence between isomorphism classes
of pairs $(A,\p)$ with $\deg \p=e$ as above and normal
algebraic affine surfaces $V$ equipped with an effective
$G_e$-action up to equivariant isomorphism.
\eprop

Thus to describe normal affine surfaces with a $G_e$-action up
to equivariant isomorphism we classify in this section
all elliptic, parabolic and hyperbolic pairs $(A,\p)$ with
$e=\deg
\p$. Our main results are the structure theorems \ref{linact},
\ref{strthm},
\ref{TN},
\ref{main} and Corollary \ref{maincor}. It also turns out
that in many cases the isomorphism class of a pair $(A,\p)$
depends only on the isomorphism class of the graded algebra
$A$, see Proposition \ref{unique}.

\subsection{Elliptic case}
Let $(A,\p)$ be an elliptic pair. It is shown in \cite[Lemmas 2.6
and 2.16]{FlZa2} that $A\cong \C[X,Y]^{\Z_d}$, where  $\C[X,Y]$
is graded via $\deg\,X=d_X>0$, $\deg\,Y=d_Y>0$, and where
$G:=\Z_d$ acts homogeneously on $\C[X,Y]$. In
particular $V=\Spec A$ is a toric surface. Moreover,
$\p$ extends to a homogeneous locally nilpotent derivation also
denoted by
$\p:\C[X,Y]\to \C[X,Y]$, and the actions of $\p$ and $G$ on
$\C[X,Y]$ commute (see Lemma \ref{obser}(a)).

\bthm\label{linact}
If $(A,\p)$ is an elliptic pair then, after an
appropriate change of coordinates, we have $A=\C[X,Y]^{\Z_d}$
with $G=\Z_d=\langle\zeta\rangle$, where $\zeta$ is a primitive
$d$-th root of unity generating $G$, acting on $\C[X,Y]$ via
$$
\zeta.X=\zeta X,\qquad \zeta.Y=\zeta^e Y\,,
$$
and $\p$ extends to $\C[X,Y]$ via
$$
\p(X)=0,\qquad
\p(Y)=X^e\qquad \mbox{i.e.,}\quad \p=X^e{\p\over\p Y}\,,
$$
where $e\ge 0,\,\,\,\gcd\,(e,d)=1$.
\ethm

\bproof Since $\p$ is locally nilpotent on $\C[X,Y]$ we have
$\p(P)=0$ for an irreducible quasihomogeneous polynomial
$P\in\C[X,Y]$ with $\deg\,P>0$ (see Proposition \ref{Ren}(b)). We
can write $\p=P^s{\tilde\p}$, where ${\tilde \p}$ is again a
locally nilpotent derivation and $s$ is chosen to be maximal. The
derivation, say, $\bar \p$ of $\C[X,Y]/(P)$ induced by ${\tilde
\p}$ is then nontrivial, so by Proposition \ref{Ren}(c)
above $\C[X,Y]/(P)$ is a polynomial ring in one variable. Since
$P$ is quasihomogeneous, it must be linear in $X$ or in $Y$.
After a suitable quasihomogeneous change of variables we may
assume that
$P=X$ so that $\p(X)=0$ and $\ker\p=\C[X]$. Since $\p$ is
homogeneous locally nilpotent, $\p(Y)$ is a homogeneous
polynomial in $X$, i.e., $\p(Y)=aX^e$ with $a\in \C^*$ and $e\ge
0$ (cf. the proof of Lemma 2.16 in \cite{FlZa2}). Replacing $Y$
by $Y/a$ we may suppose that $a=1$.

Since $\p$ commutes with the action of $G$, for any $g\in G$ we
have $\p(g.X)=g.\p(X)=0$, and so $g.X=\alpha(g)X$ for some
character $\alpha : G \to S^1$. It was shown in the proof of
\cite[Lemma 2.16]{FlZa2} that $\alpha$ is necessarily injective.
Thus we can identify $G$ with the cyclic group
$\alpha(G)=\langle\zeta\rangle\cong\Z_d$ for a certain primitive
$d$-th root of unity $\zeta$, where $\zeta .X=\zeta X$. We write
now $\zeta.Y=\alpha Y+\beta X^{\sigma}$, where $d_X=\sigma d_Y$ in
the case that $\beta\neq 0$. Since $\p(\zeta.Y)=\zeta.\p(Y)$ we
obtain
$$\alpha X^e=\zeta.X^e=\zeta^e X^e\,,$$
and therefore $\alpha=\zeta^e$. If $\gcd\,(d,e)\neq 1$
then $d'e\equiv 0\mod d$ for some $d'<d$, and so $\zeta^{d'}\neq
1$ acts as a pseudo-reflection on $\C[X,Y]$, which is excluded by
our assumption that $G$ is small. Hence $\gcd\,(d,e)=1$.

Finally, if $\zeta^e=\zeta^{\sigma}$ then $\zeta$ when considered
as an operator on $\C Y+\C X^{\sigma}$ has infinite order, which
is impossible. Hence $\beta=0$ in this case. If $\zeta^e\neq
\zeta^{\sigma}$ then replacing $Y$ by $Y':=Y+{\beta\over
\zeta^e-\zeta^{\sigma}} X^{\sigma}$ we can achieve that
$\zeta.Y'=\zeta^e Y'$, proving the theorem. \eproof

\subsection{Technical lemmas}

\bnota\label{grdm} Until the end of this section we let $(A,\p)$
be a parabolic or hyperbolic pair as in Definition \ref{ellpa}.
Thus $\p$ is a homogeneous locally nilpotent derivation on
$A=A_+\oplus A_0\oplus A_-$ corresponding to a $\C_+$-action,
$C:=\Spec A_0$ is a smooth curve and $A_+:=\bigoplus_{i>0}
A_i\neq 0$. We assume as before that the $\C^*$-action is
effective so that $A_1\neq 0$, and also $A_{-1}\neq 0$ as soon as
$A_{-}:=\bigoplus_{i<0} A_i\neq 0$. We let $d=d(A_{\ge 0})$ be
the minimal positive integer such that $A_{d+n}=A_dA_n$ for every
$n\ge 0$ (see \cite[3.6 and Lemma 3.5]{FlZa1}).
   \enota

\blem\label{A} If $\p|A_0\neq 0$ then $A_0=\C[t]$ for a certain
$t\in A_0$. Consequently for every $m\in M,\,\,\,$ $A_m$ is a
       free $A_0$-module of rank
$1$.\elem

\bproof The morphism $\pi:\Spec A \to C=\Spec A_0$ induced by the
inclusion $A_0\hookrightarrow A$ coincides with the orbit map
onto the algebraic quotient $V\quot\C^*$, hence its general fiber is
an orbit closure of the $\C^*$-action on $V=\Spec A$ associated
to the given grading. Since $\p|A_0\neq 0$ the general orbits of
the $\C_+$-action $\varphi_{\p}$ on $V$ belonging to $\p$ are not
contained in the fibers of $\pi$, and so map dominantly onto
$\Spec A_0$. These orbits being isomorphic to $\A^1_\C$, $A_0$ is
a subring of a polynomial ring $\C[T]$. It is easily seen that
$A_0$ is a normal ring, hence $A_0=\C[t]$ for some $t\in A_0$, as
stated. Now the second statement follows from \cite[Lemma
1.3(b)]{FlZa1}.\eproof

For later use we consider in the next lemma more
generally a non-homogeneous derivation, but with homogeneous
components of only nonnegative degrees.

\blem\label{kernel}
\footnote{Cf. Lemma \ref{quo}.}
Let $\delta=\sum_{i=k}^l \p_i$ be a nonzero
locally nilpotent derivation on $A$ decomposed into
homogeneous components with $l\ge k\ge 0$. If $d:=d(A_{\ge
0})$ and $v\in A_d$ generates $A_d$ as an $A_0$-module, then
$\ker\delta=\C[v, v^{-1}]\cap A$. In particular, $\p|A_0\ne 0$.
\elem

\bproof
Note first that $\delta$
stabilizes the subring $A_{\ge 0}$. Since by definition of $d$
we have $A_{n+d}=A_nA_d=vA_n$, it stabilizes as well the
principal ideal $vA_{\ge 0}$ of $A_{\ge 0}$. Thus by Corollary
\ref{norm}(b) $\delta(v)=0$ and so $\C[v,
v^{-1}]\cap A\subseteq \ker\delta$. To deduce the other
inclusion it is sufficient to show that $\C[v,
v^{-1}]\cap A$ is integrally closed in $A$ (see Proposition
\ref{Ren}(b)). The normalization of $\C[v,
v^{-1}]\cap A$ in $A$ is again graded and normal and so is
equal to $\C[w,w^{-1}]\cap A $ for some homogeneous element
$w\in A$ of positive degree $d'$. Thus
$v=cw^k$ for some $k\ge 0$ and $c\in\C$, and so $d=d'k$. It
follows that
$A_{n+d}=vA_n=wA_{n+(k-1)d'}$ for all $n\ge 0$. By
definition of
$d$, this is only possible in the case
$d=d'$, which proves that $\C[v,
v^{-1}]\cap A=\ker \delta$.
\eproof

This lemma has the following important consequence. Although
it also follows from the classification theorems \ref{TN}
and \ref{main} we give here an independent proof.

\bprop\label{unique}
Let $A$ be a parabolic or hyperbolic algebra as above and let
$\p,\p'$ be nonzero homogeneous locally nilpotent derivations
on $A$ of the same degree $e$. In the parabolic case assume
further that $e\ge 0$. Then $\p$ and $\p'$ are proportional,
i.e.\
$\p'=c\p$ for some $c\in \C^*$. In particular, the pairs
$(A,\p)$ and $(A,\p')$ are isomorphic.
\eprop

\bproof
If $A$ is hyperbolic we may reverse the grading, so in
both cases we may suppose that $e\ge 0$. By the preceding lemma
$\ker\p=\ker\p'=\C[v,v^{-1}]\cap A$, where $v$ is as above.
Thus $v:V:=\Spec A\to \Gamma:=\Spec (\C[v,v^{-1}]\cap A)$ is an
affine ruling (see also Lemma \ref{quo}), and the vector fields
$\p$ and
$\p'$ are both tangent to the fibres of $v$. Hence $\p'=c\p$
for some $c\in \Frac(A)$ of degree 0, and because of Proposition
\ref{Ren}(b) we have
$c\in\ker \p$. By Lemma \ref{kernel} this implies that $c\in\C$,
proving the first assertion.

To deduce the second one, we write $c=\lambda^e$ with
$\lambda\in \C^*$. The $\C^*$-action then gives an isomorphism
$\lambda.-:A\to A$ with $\lambda.\p'=c\p'=\p$, as required.
\eproof

\blem\label{B}
If $\deg\p\ge 0$ and if $\p(u)=0$ for some
nonzero element $u\in A_{1}$, then $A_{\ge 0}\cong \C[t,u]$
with $\deg t=0$, and $\p |A_{\ge 0}=x\p/\p t$ for some
homogeneous $x\in A_{\deg \p}$. \elem

\bproof
First we note that $\p(x)\ne 0$ for all $x\in
A_0\backslash\C$ by Lemma
\ref{kernel}. Applying Lemma \ref{A} we
see that
$A_0=\C[t]$ for some
$t\in A_0$ and, moreover, for every $k>0$ the $A_0$-module
$A_k$ is freely generated by some element $e_k\in A_k$.
Therefore
$u^k=p_ke_k$ for a certain $p_k\in A_0$. Since $u^k\in \ker \p$
and $\ker \p$ is factorially closed, $\p (p_k) = \p (e_k) =0$.
Hence $p_k\in \C$ for all $k>0$, and so $A_{\ge 0}=\C[t,u]$.
Since $\p(u)=0$ we have $\p|A_{\ge 0}=x\p/\p t$, where $x=\p (t)
\in A_{\deg \p}$, as required.\eproof

\blem\label{DD}
If $\deg \p =:e\ge 0$ then there is an
isomorphism of graded $\C$-algebras
$A_{\ge 0}\cong \C[s,u']^{\Z_d}$ with $s^d=t$ and
$u^{\prime d}=v$,  where the polynomial ring $B:=\C[s,u']$
is graded via
$\deg s=0$, $\deg u'=1$ and the cyclic group
$\Z_d=\langle\xi\rangle$ acts on $B$ via
$$
\xi . s=\xi^e s,\qquad \xi . u' = \xi u'.
$$
Moreover $\gcd (e,d)=1$, and $\p$ is the restriction to $A_{\ge
0}$ of the derivation
$$
\p=u^{\prime e} {\p\over \p s}\,.
$$
\elem

\bproof
We may suppose that
$A=A_{\ge 0}$, and we let $B$ be the normalization of $A$ in
the field of fractions of $A[u']$, where $u':={\root d \of v}$.
In view of the minimality of $d$ the assumptions of
Lemma \ref{cyex} are fulfilled. Hence the
group $\Z_d\cong
\langle\xi\rangle$ acts on $B$ via $\xi | A =$id,
$\xi.u'=\xi u'$, so that $A=(B)^{\Z_d}$, and
$\p$ extends to a locally nilpotent derivation
(also denoted $\p$) on $B$ of degree $e$.
As $\p(u')=0$ and $\deg
u'=1$ we can
apply Lemma \ref{B} to obtain that $B\cong \C[s,u']$
for some $s\in B_0$, and $\p=x\p/\p s$ for
a certain homogeneous element $x=p(s)u^{\prime e}\in
\left(\C[s,u']\right)_e$.  Since $\p=p(s)u^{\prime e}
\p/\p s$ is locally nilpotent we have $p\in\C^*$. Hence we may
assume that $x=u^{\prime e}$.

The action of $\Z_d$ on $\Spec B_0 = \Spec\C[s]$ has a fixed
point which we may suppose to be given by $s=0$. Thus
$\xi.s=\xi^ks$ for some $k\in \Z$. Since $\p$ commutes with the
action of $\Z_d$ (see Lemma \ref{obser}) we have
$$\xi^e.u^{\prime e}=\xi.\p (s)=
\p (\xi.s)=\xi^ku^{\prime e}\,,$$ i.e., we may assume that $k=e$.

Since $A_1=(B_1)^{\Z_d}\neq 0$ there exists a non-zero element
$f=q(s)u'\in A_1$, where $q(s)=\sum_{m=0}^n q_ms^m\in\C[s]$.
Since $f$ is invariant under $\xi$ we obtain
$$
\xi.f=q(\xi^e s)\xi u'=q(s)u'=f\,,
$$
i.e., $\xi^{me+1}=1$ as soon as $q_m\neq 0$. Thus $me+1\equiv
0\!\mod d$ and so $\gcd(e,d)=1$. Finally, by Lemma \ref{A},
$s^d\in \C[s]^{\Z_d}=\C[t]=A_0$ generates
$A_0$. After rescaling we may suppose that $s^d=t$ as claimed.
\eproof

\brem\label{remDD} In the situation of Lemma \ref{DD}
$\Frac \left(A[{\root d \of t}]\right)=
\Frac \left(A[{\root d \of v}]\right)=\C(s,u')$.
\erem

\subsection{Parabolic case}
We are now in position to exhibit
the structure of
$(A,\,\p)$ in the case of a positive
grading with $\dim A_0=1$. We
distinguish the following cases.

\bdefi\label{PP}
A parabolic pair $(A,\p)$ as in Definition \ref{ellpa}
will be called
{\it vertical} or of {\em
fiber type} if $\p|A_0= 0$, and of
{\em horizontal type} if
$\p|A_0\ne 0$.
\edefi

Two isomorphic pairs $(A,\p)$ and
$(A',\p')$ have the same
{\it numerical invariants} $(d,e)$,
where $e:=\deg\p$ and
$d:=d(A)$ is as in \ref{grdm}
(see also \cite[3.6]{FlZa1}).
In Theorem \ref{TN} below we
show the converse, namely, that two
parabolic pairs of horizontal
type with the same numerical invariants
are isomorphic.

A parabolic pair
is of fiber type if and
only if the general orbits of the
corresponding $\C_+$-action on
$V=\Spec A$ coincide with the general
fibers of the morphism $\pi
: V \to C:=\Spec A_0$ or, equivalently,
if the vector field $\p$
on $V$ is tangent to the fibers of $\pi$.
In contrast, if the
pair is of horizontal type then the fibers
of the $\C_+$-action
map surjectively onto the base curve $C$
and so, $C\cong \A^1_\C$
or, equivalently, $A_0\cong \C[t]$
(see Lemma \ref{A}).

\smallskip

We start with the case of parabolic pairs
of fiber type.

\bthm\label{strthm} If $(A,\p)$ is a parabolic pair  of fiber
type, then $\p$ has degree $-1$. Furthermore, if we represent
$A$ via the DPD construction as
$$
A\cong A_0[D]=\bigoplus_{n\ge 0} H^0(C,
{\cO}_C(\lfloor nD
\rfloor))\cdot u^n\subseteq \Frac(A_0)[u]
$$
with a $\Q$-divisor $D$ on $C=\Spec A_0$
then $\p$ extends
to $\Frac(A_0)[u]$ as
$\p=\varphi
\frac{\p}{\p u}$, where $\varphi=\p u$
belongs to $H^0(C,
{\cO}_C(\lfloor -D \rfloor))$. Vice versa,
any $\varphi\in H^0(C,
{\cO}_C(\lfloor -D \rfloor))$ gives rise
to a homogeneous locally
nilpotent derivation $\p=\varphi \frac{\p}{\p u}$
on $A$
of degree $-1$. \ethm

\bproof
The case $\deg \p \ge 0$ is impossible by Lemma
\ref{kernel}. If
$\deg\p<0$ then $A_0\subseteq \ker\p$, and since $A_0$ is
integrally closed in $A$ we have even equality (see
Proposition \ref{Ren}(b)). If $\deg \p<-1$ then any element
in $A_1$ would be in $\ker\p$, which is a contradiction.
It follows that $\deg
\p=-1$.

If $\varphi$ is a section in $H^0(C, {\cO}_C(\lfloor -D
\rfloor))$ then the derivation $\p=\varphi \frac{\p}{\p u}$ of
$\Frac(A_0)[u]$ stabilizes $A$. Indeed, for $f\in H^0(C,
{\cO}_C(\lfloor nD \rfloor))$ we have $\varphi f \in H^0(C,
{\cO}_C(\lfloor (n-1)D \rfloor))$ and so $\p(fu^n)=n\varphi f
u^{n-1}\in A_{n-1}$. Conversely, if $\p$ is a $A_0$-linear
derivation of $A$ then it extends to $\Frac(A_0)[u]$, and so is
of type $\p=\varphi \frac{\p}{\p u}$ for some $\varphi\in
\Frac(A_0)$. If $d\in \N$ is such that $dD$ is integral then
multiplication by $\varphi$ gives a map
$$
H^0(C, {\cO}_C(\lfloor dD \rfloor))
\lto H^0(C, {\cO}_C(\lfloor
(d-1)D \rfloor))\,,
$$
and hence amounts to a section in
$H^0(C, {\cO}_C(\lfloor -D
\rfloor))$. \eproof

\brems\label{adlnd} 1. Our proof shows that

\smallskip

(i) {\it $A\cong A_0[D]$ always admits a non-zero locally
nilpotent derivation of fiber type}, and

\smallskip

(ii) {\it every homogeneous locally nilpotent derivation on $A\cong
A_0[D]$ of negative degree has degree $-1$ and is of fiber type}.

\smallskip

\noindent (i) also follows from Lemma \ref{afru}, as for a
parabolic $\C^*$-surface $V=\Spec A_0[D]$ the canonical
projection $\pi : V \to C=\Spec A_0$ is an affine ruling.

We claim as well that

\smallskip

(iii) {\it The reduced fibers of the affine ruling $\pi: V\to C$
are all irreducible and isomorphic to
$\A_\C^1$}.

\smallskip

To show (iii), with the same argument as in the proof of
Proposition 3.8(b) in \cite{FlZa1} we can reduce to the case that
$A_0=\C$ (i.e., $C=\A_\C^1$) and $D=-\frac{e'}{d}[0]$, where
$0\le e' <d$ and $\gcd(e',d)=1$ (see \cite[Theorem
3.2(b)]{FlZa1}).  In this case the reduced fiber of $\pi: V\to
\A_\C^1$ over $0\in\A_\C^1$ is isomorphic to $\Spec \C[v]$ with
$v:=t^{e'}u^d$. In fact, using the presentation of $A$ as in
(\ref{tor1}) it is readily seen that the radical of $\sqrt{tA}$
is given by
$$
\sqrt{tA} \cong \bigoplus_{b\ge 0,ad-be'> 0} \C
t^au^b,\quad\mbox{and so }\quad A/\sqrt{tA}\cong \bigoplus_{b\ge
0,ad-be'= 0} \C t^au^b\cong\C[v]\,.$$

2. The multiple fibers of $\pi: V\to C$
correspond to the points in
$|\{D\}|$. More precisely, if $\{D\}=\sum_i \frac{e_i}{m_i} a_i$ with
$a_i\in C$ and $\gcd(e_i,m_i)=1$ then $\pi^*(a_i)=m_i \pi^{-1}(a_i)$
(see
\cite[Theorem 4.18]{FlZa1}).

3. Let $W=\Spec B$ be any affine surface with a non-trivial
$\C_+$-action. The coordinate ring $B$ is filtered by the kernels
$B_n:=\ker \p^n$, where $\p\in \Der\, B$ is the corresponding
locally nilpotent derivation. Consider the associated graded ring
$A:=\bigoplus_{i\ge 0} A_i$ with $A_i:=B_{i+1}/B_{i}$ and the
associated homogeneous locally nilpotent derivation $\p'\in \Der\,
A$ of degree $-1$. Then $\p' | A_0=0$, and so the normalization
of $A$ is as in Theorem \ref{strthm}. \erems

In the following example we exhibit
a particular family of
parabolic pairs of horizontal type,
and then we show in Theorem
\ref{TN} below that this family
is actually exhaustive.

\bexa\label{pap} Given coprime integers $e\ge 0$ and $d>0$ let
$e'$ be the unique integer with $0\le e'<d$ and
$ee'\equiv 1\mod d$; we note that by this condition $e'=0$ and
$d=1$ if $e=0$. Letting $A_0=\C[t]$, we consider the
$A_0$-algebra $A$ given by the DPD construction as follows:
$$
A:=A_0\left[D\right]\subseteq \Frac (A_0)[u]\qquad
\mbox{with}\qquad D=-\frac{e' }{d}[0]\in {\rm Div}\ (\A^1_\C)\,.
$$ Clearly $d=d(A)$ (see Lemma 3.5 in \cite{FlZa1}).
According to \cite[Proposition
3.8]{FlZa1} and Example \ref{tosu} above we can represent
$A$ as the ring of invariants
$$
A=A^{\prime\Z_d}\qquad\text{with}\qquad A':=\C[s,u'],\quad \deg s=0,
\quad \deg u'=1 \,,
$$
where $s^d=t$, $u'=us^{e'}$, and where $\Z_d=\langle\zeta\rangle$
acts on $A'$ via
$$
\zeta . s= \zeta s,\qquad \zeta .
u'= \zeta^{e'} u'\,.
$$
Thus as in Example \ref{tosu} $V=\Spec A\cong V_{d,e'}$ is an
affine toric surface , and because of $ee'\equiv 1\mod d$ the
derivation
\be\label{oldco} \p':=u^{\prime e} {\p \over \p s} \in \Der A'
\ee of degree $e$ is locally nilpotent and commutes with the
$\Z_d$-action. By Lemma
\ref{obser} it restricts to a locally nilpotent derivation $\p$
of $A$. \eexa

\bdefi\label{papa} We call the pair $ P_{d,e}:= (A,\p) $ as above
{\it the parabolic $(d,e)$-pair}. \edefi

Note that $P_{d,e}$ is
of horizontal type. Moreover, two parabolic pairs $P_{d,e}$ and
$P_{\tilde d,\tilde e}$ are isomorphic if and only if $d=\tilde
d$ and $e=\tilde e$ (cf. \cite[Corollary 3.4]{FlZa1}).
In the next result we classify all
parabolic pairs of horizontal type.

\bthm \label{TN} Every parabolic pair
$(A,\p)$ of horizontal type
is isomorphic to the parabolic $(d,e)$-pair
$P_{d,e}$ with
$e:=\deg \p$ and $d:=d(A)$. \ethm

\bproof We recall (see \cite[Remark 2.5]{FlZa1})
that for $e,e'>0$
and $ee'\equiv 1 \mod d$, the $\Z_d$-actions
$G_{d,e}'$ and
$G_{d,e'}$ on $\A_\C^2=\Spec \C[s,u']$ with
$$
G_{d,e}':\quad \xi . (s,u')=(\xi^e s, \xi
u')\,\qquad\text{and}\qquad G_{d,e'}:\quad
\zeta . (s,u')=(\zeta s, \zeta^{e'} u')\,,
$$
where $\xi,\,\,\zeta$ are primitive $d$-th roots of unity with
$\xi=\zeta^{e'}$, have the same orbits, hence also the same rings
of invariants. Now Lemma \ref{DD} shows that $(A,\p)$ is
isomorphic to $P_{d,e}$. This proves the
result. \eproof

\bexa\label{deg0} If $A$ is parabolic and admits a nonzero
homogeneous locally nilpotent derivation $\p$ of degree 0 then
$A\cong\C[t, u]$ and $\p=\frac{\p}{\p t}$. In fact, by the
classification above $(A,\p)$ is the pair $P_{1,0}$ i.e., $e'=0$,
$d=1$, $s=t$ and $u'=u$ in Example \ref{pap}. \eexa

\brems\label{mulpar} 1. We note that the derivation $\p$ of in
Example \ref{pap}
  naturally extends to $\Frac (A_0)[u,u^{-1}]$ giving the
derivation \be\label{newco} \p=d\cdot t^{k+1}u^e {\p\over\p t}-
e'\cdot t^{k}u^{e+1} {\p\over\p u} =t^ku^e\left ( d\cdot t
{\p\over\p t}-e'\cdot u {\p\over\p u}\right )\,, \ee where
$ee'-1=kd$. Indeed, from $t=s^d$ and $u=u's^{-e'}$ we obtain
$$
\p (t)=d\cdot s^{d-1} u^{\prime e} =d\cdot
t^{k+1}u^e\quad\text{and}\quad \p (u)=-e'\cdot s^{-e'-1}u^{\prime
e+1}=-e'\cdot t^ku^{e+1}\,.
$$

2. By virtue of Lemma \ref{kernel}, $\ker \p =\C[v]$. Hence
$v: V\to \A_\C^1$ is the orbit map of the $\C_+$-action $e^{t\p}$
on $V$ generated by $\p$. As $v$  is homogeneous of degree
$d=d(A)>0$, the $\C^*$-action on $V$ acts non-trivially on this
affine ruling and on its base. Therefore $v$ can have at most one
degenerate fiber $v^{-1}(0)$, which is the fixed point curve
$C_+\cong\A_\C^1$ of the $\C^*$-action. Moreover, $\div (v)=dC_+$
(see \cite[Remark 3.7]{FlZa1}). \erems

\bcor\label{TNC}
A normal affine surface $V=\Spec A$, where
$A=A_0[D]$, admits a horizontal  $\C_+$-action if and only if
$A_0\cong \C[t]$ and the fractional part $\{D\}$ of the
$\Q$-divisor $D$ on $C\cong  \A_\C^1$ is supported at one point
or is zero.
\ecor

\bproof
This follows immediately from Theorem \ref{TN}; note that in
the case $A_0\cong \C[t]$  we have $A=A_0[D]\cong A_0[\{D\}]$,
see \cite[Corollary 3.4 and Proposition 3.8]{FlZa1}.

Let us provide for the `only if'-part an independent geometric
argument. For this consider
more generally a morphism
$\pi:V\to C$ of a normal affine surface $V$ onto a smooth affine
curve $C$ with only irreducible fibers. We claim that {\it if
there  exists an affine ruling $v:V\to \Gamma$ different from
$\pi$ then $C \cong\A_\C^1$ and $\pi$ has at most one multiple
fiber}. Clearly, this claim implies our assertion (see Remark
\ref{adlnd}.2). To show the claim, we let $G\cong \A_\C^1$ be a
general fiber of $v$, and we assume on the contrary that $\pi$
has at least two fibers $F_i$ of multiplicity $m_i\ge 2$,
$i=0,1$. As $\pi\vert G : G \to C$ is dominant it follows that $C
\cong\A_\C^1$, and so $\pi\vert G : G \to C$ can be viewed as a
non-constant polynomial $\nu \in\C[t]$. We also may assume that
$F_0=\pi^{-1}(0)$ and $F_1=\pi^{-1}(1)$. As $G$ is a general
fiber of $v$ it meets $F_i$ at smooth points of $V$ only, with
the intersection multiplicities in $G \cdot  \pi^{*}(i)$ being a
multiple of $m_i$ ($i=0,1$). Thus $m_0$, $m_1$, divides the
multiplicity of any root of the polynomial $\nu$,
$\nu-1$, respectively. Hence $\nu=\lambda^{m_0}=\mu^{m_1}+1$
for some non-constant polynomials $\lambda,\mu\in\C[t]$. The
pair $(\lambda,\mu)$ defines a dominant map $\A_\C^1\to
\Gamma_{m_0,m_1}$, where $\Gamma_{m_0,m_1}$ is the smooth plane
affine curve $x^{m_0}-y^{m_1}=1$. But the existence of such a map
contradicts the Riemann-Hurwitz formula, which proves our claim.
\eproof

\subsection{Hyperbolic case}
In this subsection we assume that $A$
is hyperbolic, so that $A_\pm\neq
0$. If $\p$ is a homogeneous locally nilpotent
derivation on $A$ of degree $e$ with
$e<0$ then by reversing the grading of
$A$ we obtain a
derivation of positive degree.
Thus it is sufficient to classify the
hyperbolic pairs $(A,\p)$ as in Definition \ref{ellpa}.

\blem\label{stb}
If $(A,\p)$ is a hyperbolic pair then $\p$
stabilizes $A_{\ge 0}\subseteq A$,
and $(A_{\ge 0},\p)$ is
a parabolic pair of
horizontal type.
\elem

\bproof It follows immediately from the definitions that $(A_{\ge
0},\p)$ is a parabolic pair. If it were of fiber type then the
orbits of the corresponding $\C_+$-action on $V=\Spec A$ would be
the fibers of $\pi:V\to C=\Spec A_0$. As the general fiber of
$\pi$ is $\C^*$, this leads to a contradiction.
\eproof

Thus by Theorem \ref{TN} $(A_{\ge 0},\p)$ is isomorphic to the
$(d,e)$-pair $P_{d,e}$, where $e=\deg \p$ and $d=d(A_{\ge
0})=d(\p)$ (see \ref{grdm} and Lemma \ref{DD}). In particular,
$A_0=\C[t]$ and $A_{\ge 0}= A_0[-\frac{e'}{d}[0]]\subseteq
\Frac(A_0)[u]$, where $0$ is the origin in $\A_\C^1=\Spec A_0$
(see Example \ref{pap}). Moreover $\p$ is given as in
(\ref{oldco}) or, alternatively, as in (\ref{newco}). The
following lemma is crucial in our classification.

\blem \label{twoco} Let $D_+,\ D_-$ be $\Q$-divisors on $C:=\Spec
A_0$ with $A_0=\C[t]$ satisfying $D_++D_-\le 0$, where
$D_+:=-\frac{e'}{d}[0]$ with $0\le e'<d$ and $\gcd(e',d)=1$. The
derivation $\p:A_0[D_+]\to A_0[D_+]$ of degree $e\ge 0$ as in
(\ref{newco}) extends to
$$
A=A_0[D_+, D_-] \subseteq\Frac(A_0)[u,u^{-1}]
$$
if and only if the following two conditions are
satisfied.
\begin{enumerate}
\item[(i)] If $D_-(0)\ne \frac{e'}{d}$ then
$\frac{ee'-1}{d}-eD_-(0)\ge 0$ i.e.,
$-e(D_+(0)+D_-(0))\ge 1/d$.
\item[(ii)] If $a\in\A_\C^1$ with $a\ne 0$ and
$D_-(a)\ne 0$ then
$-1-eD_-(a)\ge 0$.
\end{enumerate}
\elem

\bproof Note that $\p$ extends in a unique way
to a derivation of
$\Frac(A_0)[u,u^{-1}]$ also denoted $\p$.
We must show that $\p$
stabilizes $A$ if and only if (i) and (ii)
are satisfied.

Let us first treat the case $d=1$ so that $e'=0$ and $D_+=0$.
Then (i) and (ii) can be reduced to the condition \be\label{miq}
-1-eD_-(a)\le 0 \qquad\forall a\in \A^1_\C\,.\ee Moreover, $k=-1$
and so according to (\ref{newco}) $\p=u^{-e}\frac{\p}{\p t}$ acts
on a homogeneous element $f(t)u^{-m}\in\C(t)u^{-m}$ by
\be\label{via} \p(f(t)u^{-m})= f'(t)u^{e-m} \,.\ee  Thus $\p$
stabilizes $A$ if and only if $f(t)\in A_{-m}u^{-m}$ ($m\ge 0$)
implies $f'(t)u^{e-m}\in A_{e-m}$ or, equivalently, \be
\label{eq0} \div\, f+mD_-\ge 0\Rightarrow\left\{
\begin{array}{ll}
\div\,f'+(m-e)D_-\ge 0&\mbox{if }m-e\ge 0\\
\div\,f'+(e-m)D_+\ge 0&\mbox{if }m-e\le 0\,.
\end{array}\right.
\ee
If (\ref{miq}) is satisfied then for any
$a\in\A^1_\C$
$$
\div_a f'+(m-e)D_-(a)\ge \div_a f+mD_-(a)-1-eD_-(a)
\ge \div_a
f+mD_-(a)\,,
$$
where $\div_a(\ldots)$ denotes the order at
$a$. Thus (\ref{eq0}) is satisfied if $m-e\ge 0$,
and since
$D_+(a)=0$ and $-D_-(a)\ge 0$, it also follows
     for $m-e\le 0$.

Conversely, assume that $\p$ stabilizes $A$.
Consider $m>e$ such
that the divisor $mD_-$ is integral.
For $a\in \A^1_\C$ with
$D_-(a)\ne 0$ we let $s:=-mD_-(a)$;
thus $s\ge 0$.
Consider a polynomial $Q$
without zero at $a$ such that
$$
(t-a)^sQ u^{-m}\in A_{-m}.
$$
By assumption
$\p((t-a)^sQu^{-m})=
(s(t-a)^{s-1}Q+(t-a)^sQ')u^{e-m}\in A_{-m+e}$
and so
$$
\div_a(s(t-a)^{s-1}Q+(t-a)^sQ')+
(m-e)D_-(a)\ge 0\,.
$$
The term on the left is equal to $s-1$,
hence we obtain
$$
s-1+ (m-e) D_-(a) = -1-eD_-(a)\ge 0\,,
$$
as required in (\ref{miq}).

In  case $d\ge 2$ we consider the normalization $A'$ of $A$ in
$\Frac(A)[\sqrt[d]{v}]$ as in Lemma \ref{DD}, and we let
$p:\A^1_\C\cong \Spec A_0'\to \A^1_\C \cong\Spec A_0$,
$s\longmapsto s^d$, be the covering induced by the inclusion
$A_0\subseteq A_0'$. By {\it loc.cit.}\ $A'_{\ge 0}\cong
\C[s,u']$ with $s^d=t$ and $u^{\prime d}=v\in \ker \p$, $\deg
u'=1$, and $\p$ extends to the derivation $\p'=u^{\prime
e}\frac{\p}{\p s}$ on $A'_{\ge 0}$ and as well on $\Frac(A')$. If
$\p$ stabilizes $A$ then $\p'$ stabilizes $A'$ (see Lemma
\ref{cyex}). Moreover, $v$ can be written as $v=t^{e'} u^d$ (see
the proof of Theorem 4.15 in \cite{FlZa1}). So, by
\cite[Proposition 4.12]{FlZa1},
$$
A'\cong A_0'[D_+',D_-']\subseteq \Frac(A_0')
[u',u^{\prime -1}]\,,
$$
where $D'_+=0$ and $D'_-=p^*(D_++D_-)$.
Using the first part of
the proof we get that $D'_-(a')<0$ implies
$-1-eD'_-(a')\ge 0$.
If $a=p(a')\neq 0$ then $D_-(a)=D_-(a')$,
hence (ii) follows.
Similarly, if $p(a')=0$ then
$D_-'(a')=-e'+dD_-(0)$ and (i)
follows.

Conversely, assume that (i) and (ii) are satisfied. Reversing the
reasoning above we obtain that $-1-D'_-(a')\ge 0$ if $D'_-(a')\ne
0$. Thus by the first part $\p'$ stabilizes $A'$. Taking
invariants $\p$ stabilizes $A=(A')^{\Z_d}$, as desired.
\eproof

Summarizing we state now our main classification result for
hyperbolic pairs.

\bthm\label{main} If $(A,\p)$ is a hyperbolic pair with
$d:=d(A_{\ge 0})$ and $e:=\deg \p$, then $A_0\cong \C[t]$ and
$A\cong A_0[D_+,D_-]$ for two $\Q$-divisors $D_+$, $D_-$ on
$\A^1_\C$ with $D_++D_-\le 0$, where the following conditions are
satisfied:

\begin{enumerate}
\item[(i)] $D_+=-\frac{e'}{d}[0]$ with $0\le e'<d$
and $ee'\equiv 1\mod d$.
\item[(ii)] If $D_+(a)+D_-(a)\neq 0$ then $-(D_+(a)+D_-(a))^{-1}\le
\left\{\begin{array}{lll} de,\quad  a=0\\ e,\quad a\neq 0\\
\end{array}
\right.$.
\item[(iii)] $\p$ is defined by (\ref{newco})
in
Remark \ref{mulpar}.
\end{enumerate}
Conversely, given two $\Q$-divisors $D_+$ and
$D_-$ on $\A^1_\C$
with $D_++D_-\le 0$ satisfying (i) and (ii)
there exists a unique, up to a constant,
locally nilpotent derivation $\p$ of degree $e$
on
$A=A_0[D_+,D_-]$, and this derivation is  as in (iii).
In
particular, isomorphism classes of hyperbolic
pairs are in
1-1-correspondence to pairs $(P_{d,e}, D_-)$,
where $D_-$ is a
$\Q$-divisor on $\A^1_\C$ satisfying (ii).
\ethm

\bproof By Theorem \ref{TN}, $(A_{\ge 0},\p)$ is isomorphic to the
parabolic pair $P_{d,e}$. In particular, (i) and  (iii) are
satisfied. By Lemma \ref{twoco} also (ii) holds, proving the
theorem. \eproof

\bcor\label{crucol} A two-dimensional  normal graded $\C$-algebra
$A=\bigoplus_{m\in \Z} A_m$ with $A_\pm\neq 0$ admits a
homogeneous locally nilpotent derivation $\p$ of positive degree
if and only if $A_0\cong \C[t]$ and $A\cong
A_0[D_+,D_-]$, where the fractional part
$\{D_+\}$ is supported at one point or is zero.
\ecor

In order to study more closely the structure of the affine ruling
which corresponds to the $\C_+$-action with generator $\p$ as
above, we need a simple lemma. We let $A=A_0[D_+,D_-]\subseteq
\Frac (A_0)[u,u^{-1}]$ be a normal graded $\C$-algebra, and we
consider the associated $\C^*$-fibration $\pi:V=\Spec A\to
C:=\Spec A_0$ over the curve $C$. It was shown in \cite[Theorem
4.18]{FlZa1} that the fiber over a point $a\in C$ with
$D_+(a)+D_-(a)<0$ consists of two $\C^*$-orbit closures
$\bO^\pm_{a}$. Moreover, if $D_+(a)=-\frac{e_+}{m_+}$,
$D_-(a)=\frac{e_-}{m_-}$, where $m_+>0$, $m_-<0$ and $\gcd(e_\pm,
m_\pm)=1$, then \be\label{*} \qquad
\pi^*(a)=m_+[\bO^+_{a}]-m_-[\bO^-_{a}]\qquad\mbox{and}\qquad
\div\, u = -e_+[\bO^+_{a}]+e_-[\bO^-_{a}] +\ldots\,,\ee where the
terms in dots correspond to points in $|D_+|\cup |D_-|$ different
from $a$. Letting $v_\pm\in A_{m_\pm}$ be an element with
$A_{m_\pm}=v_\pm A_0$ near $a$, we have the following observation.

\blem\label{trmap} \begin{enumerate} \item[(a)] The orbit closures
$\bO^\pm_{a}\cong \Spec \C[v_\pm]$ are smooth affine lines.
\item[(b)] $\div
(v_+)=\Delta(a)[\bO^-_{a}]$ and $\div
(v_-)=\Delta(a)[\bO^+_{a}]$, where $\Delta(a):=m_+e_--m_-e_+$.
\end{enumerate}
\elem

\bproof With the same argument as in the proof  of Proposition
3.8(b) in \cite{FlZa1} we can  reduce to the case that $A_0=\C[t]$
and $\vert D_+\vert\cup\vert D_-\vert$ is the point $a=0\in
C=\A^1_\C$. We may also suppose that $D_+(0)+D_-(0)<0$. Recall
(see the proof of Theorem 4.15 in \cite{FlZa1}) that
$v_\pm=t^{e_\pm} u^{m_\pm}$ up to a constant in $\C^*$.

(a) The ideal of $\bO_a^+$ coincides with the radical
$\sqrt{v_-A}$, see the proof of Theorem 4.18 in \cite{FlZa1}. Thus
it suffices to show that
$$\sqrt{v_-A}=A_-\oplus
tA\oplus\bigoplus_{m_+\nmid \beta} A_{\beta}\,.$$ As $v_+\not\in
\sqrt{v_-A}$ we have the inclusion '$\subseteq$'. To deduce
    '$\supseteq$' we note first that
$A_-\oplus tA\subseteq \sqrt{v_-A}$. Suppose that
$t^{\alpha}u^{\beta}\in A$, where $\beta>0$  and $m_+\nmid
\beta$, and let us show that
$t^{\alpha}u^{\beta}\in\sqrt{v_-A}$. For this we need to prove that
$v_-=t^{e_-}u^{m_-}$ divides $t^{n\alpha}u^{n\beta}$ in $A$ for
$n\gg 0$ or, equivalently, that $t^{n\alpha-e_-}u^{n\beta-m_-}\in
A_{n\beta-m_-}$. This amounts to
\be\label{radi}\qquad
(n\alpha-e_-)+(n\beta-m_-)D_+(0)\ge 0\,\,\,\Leftrightarrow\,\,\,
n(\alpha+\beta D_+(0))\ge e_-+m_-D_+(0)\,.\ee
Because of our assumptions $t^{\alpha}u^{\beta}\in A$ and
$m_+\nmid\beta$ we have $\alpha+\beta D_+(0) \ge0$ and $\beta
D_+(0)\not\in\Z$, so $\alpha+\beta D_+(0)> 0$. Hence (\ref{radi})
is satisfied for $n\gg 0$, as required.

(b) follows from (\ref{*}) by virtue of the equalities
$$\div \,v_+ = e_+\div \,t +m_+\div \,u\qquad\mbox{and}\qquad
\div \,v_- = e_-\div \,t +m_-\div\, u\,.$$ \eproof

We consider below a hyperbolic pair $(A,\p)$ as in Theorem
\ref{main}, and we let $v\in A_d$ be a generator of $A_d$ over
$A_0=\C[t]$ (cf. Lemma \ref{kernel}). Then $v: V=\Spec A\to
\Gamma=\A^1_\C$ provides an affine ruling which is the quotient
map of the $\C_+$-action on $V$ induced by $\p$. In the next
proposition we describe the multiplicities which occur in the
degenerate fibers of this affine ruling (cf. Remark
\ref{mulpar}.2).

\bprop\label{qmap} The fiber of the affine ruling
$v:V\to\Gamma=\A^1_\C$ over a point $p\ne 0$ is smooth,  reduced
and consists of just one $\C_+$-orbit, whereas the fiber over
$p=0$ is a disjoint union of $\C^*$-orbit closures isomorphic to
affine lines, one for each point $a\in C$ with $D_+(a)+D_-(a)<0$.
Moreover \be\label{mpc} \div (v)=d_+\sum_{a\in \A^1_\C}
m_-(a)\left(D_+(a)+D_-(a)\right) [\bO^-_{a}]\,, \ee where the
integer $m_-(a)<0$ is defined by $D_-(a)=\frac{e_-(a)}{m_-(a)}$
with $\,\,\gcd(e_-(a), m_-(a))=1$. \eprop

\bproof As $v$ is homogeneous of degree $d=d_+:=d(A_{\ge 0})$ the
affine ruling $v:V\to\A^1_\C$ is equivariant if we equip
$\A^1_\C$ with the $\C^*$-action $\lambda . t = \lambda^d t$.
This implies that for every point $p\ne 0$, the fiber of $v$ over
$p$ is smooth, reduced and consists of just one $\C_+$-orbit. By
the previous lemma, $\div (v)$ is a linear combination of the
divisors $\bO^-_{a}$, where $a$ runs through all points of
$C=\A^1_\C$ with $D_+(a)+D_-(a)<0$. We compute the multiplicities
separately in the cases where $a=0$ and $a\ne 0$.

If $a=0$ then $D_+(0)=-\frac{e_+}{d_+}$
with $e_+=e'$ and $D_-(0)=\frac{e_-}{m_-}$
with $e_-=e_-(0)$, $m_-=m_-(0)$,
so by Lemma
\ref{trmap} the coefficient of $\bO^-_{0}$ in
$\div (v)$ is $\Delta(0)=-
e_+m_-+e_-m_+=d_+m_-\left(D_+(0)+D_-(0)\right)$,
which agrees with (\ref{mpc}).

If $a\ne 0$ then $m_+(a)=1$ and so
$\Delta(a)=m_-(a)\left(D_+(a)+D_-(a)\right)$.
Letting $v_*\in A_1$
be an element generating $A_1$ over $A_0$ near $a$,
we can write $v=\varepsilon v_*^{d_+}$,
where $\varepsilon\in A_0$ is a unit near $a$
i.e., $\varepsilon(a)\ne 0$. By Lemma
\ref{trmap} $\bO^-_{a}$ occurs with multiplicity
$\Delta(a)$ in $\div (v_*)$, and so it
occurs with multiplicity $d_+\Delta(a)=
d_+m_-(a)(D_+(a)+D_-(a))$
in $\div (v)$, as required in (\ref{mpc}).
\eproof

\brem \label{nonum} We note that $\div (v)$ is the exceptional
divisor of the birational morphism $\sigma_+:V \to V_+=\Spec
A_{\ge 0}$ induced by the inclusion $A_{\ge 0}\hookrightarrow A$.
Indeed, the divisor $\div (v)=d_+C_+$ on $V_+$ is supported by
the fixed point curve $C_+\cong \A_\C^1$ of the $\C^*$-action on
$V_+$ (see Remark \ref{mulpar}). For every point $a\in C=\A^1_\C$
with $D_+(a)+D_-(a)<0$ there is a unique point $a'$ over $a$ on
$C_+$, and $\sigma_+$ is the affine modification consisting in an
equivariant blowing up of $V_+$ with center supported at all those
points
$a'\in C_+$ and deleting the proper transform of the divisor
$C_+$ (see \cite[Remark 4.20]{FlZa1}).  \erem

If $v$ is a unit in $A$ then $D_++D_-=0$,
$v:V\to\A^1_\C\backslash\{0\}$ is the quotient map,
and all
fibers of $v$ are smooth affine lines.
More precisely the following result holds.

\bcor\label{rd}
Let $(A,\p)$ be a hyperbolic pair and $d:=d(A_{\ge 0})$. If one of the
following two conditions is satisfied:
\begin{enumerate}
\item[(i)] $e:=\deg \p=0$, or
\item[(ii)] $A$ contains a unit of non-zero degree,
\end{enumerate}
then
$$
A\cong \C[z,\,v,\,v^{-1}]\qquad(\Rightarrow\qquad
V\cong\A^1_\C\times\C^*)\qquad {\rm and}\quad \p=
\p/\p z\,,
$$
where
$\deg z=-e$ and $\deg v=d$.
\ecor

\bproof In case (i) Theorem \ref{main}
(i) shows that $d=1$ and $e'=0$, so $D_+=0$, and moreover
by \ref{main}(ii) $D_+(a)+D_-(a)=0$ for all closed points $a\in
\A^1_\C$. Thus $D_+=D_-=0$ and $A=A_0[u,u^{-1}]$ for some element
$u\in A_1$. By \ref{main} (iii) and Remark \ref{mulpar}.1  $\p$
is the derivation $\p=\frac{\p}{\p t}$, which proves the result.

In case (ii), by \cite[Remark 4.5]{FlZa1}, $D_+=-D_-$, and by
Theorem \ref{main}, $D_+=-\frac{e'}{d}[0]$. Therefore, $A$ is the
semigroup algebra generated over $\C$ by all monomials $t^au^b$
with $ad-be'\ge 0$, $a,b \in\Z$ (cf. the proof of Theorem 4.15 in
\cite{FlZa1}). Choose $q\in\Z$ with $|{e\atop q}{d\atop e'}|=1$
and consider the elements
$$
v:=t^{e'}u^d\,,\quad v^{-1}\quad\mbox{and}\quad
z:=t^{-q}u^{-e}\in A\,\quad\mbox{with}\quad \deg
v=d,\,\,\deg z=-e\,,
$$
so that $u=v^{-q}z^{-e'}$ and $t=v^{e}z^{d}$. As we have noticed
above, a monomial $t^au^b=v^{ae-bq}z^{ad-be'}$ belongs to $A$ if
and only if $ad-be'\ge 0$. Thus $A= \C[v,v^{-1},z]$. The orbits
of the  $\C_+$-action on $\Spec \C[v,v^{-1},z]=\A^1_\C\times
\C^*$ given by $\p$ are necessarily contained in the fibers of
the projection to $\C^*$, and $\ker \p=\C[v,v^{-1}]$ (cf.
Lemma \ref{kernel}). Since $\p$ is homogeneous of degree
$e$, we get $\p=cv^az^b\frac{\p}{\p z}$ for suitable $c\in\C^*$,
$a\in\Z$, $b\in\N$ with $ad-be=0$. As $\p$ is also locally
nilpotent this forces $a=b=0$ and so $\p=c\frac{\p}{\p z}$.
Replacing $z$ by $z/c$, the result follows. \eproof

Next we describe explicit equations for hyperbolic pairs in the
case that $A=A_0[D_+,D_-]$ with $D_+=0$.

\bcor\label{KN} Let $(A,\p)$ be a hyperbolic pair, and suppose
that $A=A_0[D_+,D_-]$ with $D_+=0$, so that $A_{\ge 0}\cong
\C[t,u]$ with $\deg u=1$ and $\deg t=0$. If $k:=d(A_{\le 0})$ and
$e:=\deg \p\ge 0$ then $A$ is the normalization of the graded
domain
$$
B=B_{k,P}:=\C[t,u,v]/(u^{k}v-P(t))
\quad\text{with}\quad
\deg v=-k\,,
$$
where \be\label{poly} P(t)=\prod_{i=1}^s
(t-a_i)^{r_i}\in\C[t]\qquad (r_i\ge 1
\mbox{ and }a_i\neq
a_j\mbox{ for }i\neq j) \ee is a unitary
polynomial uniquely
determined by $D_-=-\div\,P/k$ and satisfying
\be\label{poly1}
\gcd(k,r_1,\ldots,r_s)=1\quad\mbox{and}\quad
e\ge k/r_i \ee
for
$i=1,\dots,s$. The derivation $\p$ is given
(and uniquely
determined) by the conditions \be\label{co}
\p(u)=0,\quad \p
(t)=u^e \quad(\Rightarrow\quad \p (v)= P'(t)u^{e-k}).
\ee
Conversely, given a polynomial $P$ as in (\ref{poly})
and
(\ref{poly1}) there is up to a constant a unique
locally
nilpotent derivation $\p$ of degree $e$ of the
normalization $A$
of $B_{k,P}$ satisfying (\ref{co}).
\ecor

\bproof As was shown in \cite[Example 4.10 and Proposition
4.11]{FlZa1}, $A$ is the normalization of the algebra $B_{k,P}$,
where $P$ is a unitary polynomial uniquely determined by
$D_-=-\div(P)/k$. Since $k$ is minimal with $kD_-$ integral, we
have $\gcd(k,r_1,\ldots,r_n)=1$. By Theorem \ref{main}(ii),
(iii) it follows that $e\ge k/r_i$ and that $\p$ has the stated
form (\ref{co}). Conversely, given $P$ the normalization $A$ of
$B_{k,P}$ is isomorphic to $A_0[D_+,D_-]$ with $D_+=0$ and
$D_-=-\div(P)/k$. If $e\ge k/r_i$ for all $i$ then the conditions
(i), (ii) in  Theorem \ref{main} are fulfilled for $A$, so there
is a locally nilpotent derivation $\p$ of $A$ satisfying
(\ref{co}), and $\p$ is uniquely determined up to a constant
factor. \eproof

\brems\label{aut} 1. Over each of the points $t=a_i\in\A_\C^1$,
the
surface  $V=\Spec A$ considered in Corollary \ref{KN} has a unique
fixed point $a_i'$ of the
$\C^*$-action. This point $a_i'\in V$ is a quotient singularity of
type
$(d_i,e_i)$, where $r_i/k=d_i/e_i'$ with $d_i,e_i'$ coprime and
$0\le e_i<d_i$, $e_i\equiv e_i'\mod d_i$. This follows from Theorem
4.15 in
\cite{FlZa1}, since $D_+(a_i)=0$ and $D_-(a_i)=r_i/k$. In particular,
the surface $V$
is smooth if and only if $r_i|k$ for all $i$ (cf. Corollary 4.16 in
\cite{FlZa1}).

2. A description of the automorphism group
$\Aut\,V_{k,P}$ for a
smooth surface $V_{k,P}:=\Spec\,B_{k,P}$,
where $B_{k,P}$ is as in
Corollary \ref{KN}, can be found in
\cite[(2.3)-(2.4)]{Be1} and
\cite[Theorem 1]{ML2}.

3. For any $e\ge k$ the derivation $\p$ described in Corollary
\ref{KN} stabilizes the ring $B$ and induces a $\C_+$-action
(actually, a $G_{e}$-action, see Lemma \ref{sdp}) on $\A_\C^3$
which leaves the surface $V_{k_,P}=\Spec\,B\subseteq \A_\C^3$
invariant. In case $e<k$, however, $\p$  does not induce a
derivation on $B$. The simplest example of such a surface
$V_{k_,P}$ is with $P=t^3$ and $k=2$, $e=1$. Here the element
$\p(v)=3t^2u^{-1}$ is not in $B$ but is integral over $B$ as
its square is equal to $9tv\in B$.

4. The $\C_+$-action associated to the derivation $\p$ in
Corollary \ref{KN} is \be\label{str} \alpha.(t,\,u,\,v)=(t+\alpha
u^e,\, u,\,u^{-k}P(t+\alpha u^e)) \,, \quad \alpha\in\C_+\, \ee
with fixed point set  $\{u=0\}$. Again, for $e\ge k$ this
$\C_+$-action extends to $\A^3_\C$.\erems

In the case $D_+=-\frac{e'}{d}[0]\ne 0$ a
suitable cyclic covering
of $V=\Spec A$ can  be described as in
Corollary
\ref{KN}. This leads to the following
alternative description of
arbitrary hyperbolic pairs $(A,\p)$.

\bcor\label{maincor} We let $(A,\p)$ be a hyperbolic pair with
invariants $d:=d(A_{\ge 0})$, $k:=d(A_{\le 0})$, $e:=\deg \p>0$.
If $A\cong A_0[D_+,D_-]$, where $D_+=-\frac{e'}{d}[0]$ and
$D_-(0)=-\frac{l}{k}$, then there exists a unitary polynomial
$Q\in\C[t]$ with $Q(0)\neq 0$ and $\div (Qt^l)=-kD_-$. Moreover if
$A'= A_{k,P}$ is the normalization of \be\quad\label{Bkp}
B_{k,P}=\C[s, u,v]/ \left( u^kv-P(s)\right)\,,
\quad\mbox{where}\quad P(s):=Q(s^d)s^{ke'+dl}\,, \ee then the
group $\Z_d=\langle\zeta\rangle$ acts on $B_{k,P}$ and also on
$A'$ via \be\label{opac} \zeta . s=\zeta s, \quad\zeta .
u=\zeta^{e'}
      u\quad\text{and}\quad
\zeta . v=v\, \ee so that $A\cong A^{\prime \Z_d}$. Furthermore,
$ee'\equiv 1 \mod d$ and, up to a constant factor, $\p$ is the
restriction of the derivation $u^e\frac{\p}{\p s}$ to $A$. \ecor

\bproof The inequality $D_+(0)+D_-(0)\le 0$ is equivalent to
$ke'+dl\ge 0$. This implies that there are unitary polynomials
$Q(t)\in\C[t]$ and $P(s)\in\C[s]$ such that $\div (Qt^l)=-kD_-$
and $P(s)=Q(s^d)s^{ke'+dl}$.

The isomorphism $A\cong A^{\prime\Z_d}$ was established in
Example 4.13 and Proposition 4.14 of \cite{FlZa1}. The derivation
$u^e\frac{\p}{\p s}$ commutes with the $\Z_d$-action (\ref{opac}),
and so restricts to a homogeneous locally nilpotent derivation
$\p'$ of degree $e$ on $A$, iff $ee'\equiv 1 \mod d$ (see Theorem
\ref{main}(i)). Thus by Corollary \ref{KN} it is equal to $\p$ up
to a constant. The rest of the proof can be left to the reader.
\eproof

\section{Applications}

\subsection{Preliminaries}
Sometimes the surfaces $V=\Spec A$
as above admit two
$\C_+$-actions with different orbit maps;
see e.g. Example \ref{tosu}. The
following example is also well known.

\bexa\label{derplmi} We let $A$ be the normalization of the ring
$B_{1,P}=\C[t,u_+,u_-]/(u_+u_--P(t))$, where $P\in \C[t]$ is a
unitary polynomial and  the grading is given by $\deg t=0$, $\deg
u_\pm=\pm 1$. By Corollary \ref{KN}, for every $e\ge 1$ there are
homogeneous derivations of degree $e$ as well as of degree $-e$
on $A$. More explicitly these are given (up to a constant factor)
by \be\label{dere} \qquad\p_{+}=u_+^e\,{\p}/{\p t}+ P'(t)
u_+^{e-1}\,{\p}/{\p u_-} \quad\mbox{and}\quad
\p_{-}=u_-^e\,{\p}/{\p t}+ P'(t) u_-^{e-1}\,{\p}/{\p u_+}\,; \ee
cf. (\ref{co}). Note that $\ker (\p_\pm)=\C[u_\pm]$, hence the
corresponding $\C_+$-actions $\varphi_{+}$ and $\varphi_{-}$
preserve the affine rulings $u_\pm : V\to \C$ of $V=\Spec A$,
respectively. These rulings are different provided that $P$ is
a non-constant polynomial.

In view of (\ref{str}) $\varphi_{+}$ is given by
$$\alpha.(t,\,u_+,\,u_-)=(t+\alpha u_+^e,\,
u_+,\,u_+^{-1}P(t+\alpha u_+^e)) \,, \quad \alpha\in\C_+\,.$$ As
$\ker (\p_{-})=\C[u_-]$ the conjugated locally nilpotent
derivation
$$\p_{\alpha}: =\alpha.\p_{-}.\alpha^{-1} \in\Der A$$
has kernel $\ker (\p_{\alpha}) =\C[u_{\alpha}]$, where
$$u_{\alpha}:=\alpha(u_-)=u_+^{-1}P(t+\alpha u_+^e)
= u_-+\sum_{j=1}^{\deg P}
P^{(j)}(t)\frac{\alpha^j}{j!}u_+^{je-1}\,.
$$
As $\alpha\in \C_+$ varies, the affine rulings $u_{\alpha}:
V\to
\A_\C^1$ also vary in a continuous family. \eexa

\bdefi\label{mli} One says that two
$\C_+$-actions on an affine
variety $V=\Spec\,A$ are {\it equivalent}
if their general orbits
are the same, or in other words, if they
define the same affine
ruling on $V$.

\edefi

If $\p$ and $\p'\in \Der\,A$
are the associated
locally nilpotent derivations then the
$\C_+$-actions are
equivalent if and only if $\ker\p=\ker\p'$,
and if and only if
$a\p=a'\p'$ for some elements $a, a'\in \ker\p$
(see \cite[Lemma
2.1]{KaML} or Proposition \ref{Ren}(b)).
Consequently, any two equivalent locally
nilpotent derivations
$\p$ and $\p'$ commute: $[\p,\,\p']=0$.

We recall \cite{KaML, Za} that the
{\it Makar-Limanov invariant} of
an affine variety $V=\Spec\,A$ is ${\rm ML}(V)={\rm
ML}(A)=\bigcap \ker\p$, where $\p$ runs
over the set of all
locally nilpotent derivations of $A$.

Certainly, a surface $V$ has a trivial Makar-Limanov invariant
${\rm ML}(V)=\C$ if and only if $V$ admits two non-equivalent
$\C_+$-actions, or two different affine rulings over affine
bases, or else two non-equivalent nonzero locally nilpotent
derivations of $A$.

A useful characterization of surfaces with a trivial
Makar-Limanov invariant is the following
result due to Gizatullin
\cite[Theorems 2 and 3]{Gi2},
Bertin \cite[Theorem 1.8]{Be2},
Bandman and Makar-Limanov  \cite{BaML3}
in the smooth case, and
to Dubouloz \cite{Du1} in the normal case.

\bthm\label{zigzag} For  a normal
affine surface $V$ non-isomorphic to $\C^*\times\C^*$,
the following conditions are equivalent.
\begin{enumerate} \item[(i)] The Makar-Limanov invariant of $V$ is trivial.
\item[(ii)] The automorphism group
$\Aut\,V$ \footnote{Which is not necessarily
an algebraic group,
see Example \ref{newgr} below.} acts on $V$
with an open orbit
$O$ such that the complement $V\backslash O$
is finite.
\item[(iii)]  $V$ admits a smooth compactification
by a {\rm zigzag} that is, by a linear chain of smooth rational
curves.
\end{enumerate}
\ethm

Thus an affine ruling $V\to\A^1_\C$ on a normal
affine surface
$V$ is unique (in other words, any two $\C_+$-actions
on $V$ are
equivalent) unless $V$ admits a smooth
compactification by a
zigzag. In the latter case there are, indeed,
at least two
different affine rulings $V\to\A^1_\C$,
hence also two
non-equivalent $\C_+$-actions on $V$.

Note that all surfaces as in Theorem \ref{zigzag} are rational and
allow a constructive description, see \cite[Proposition 3]{Gi2}
or \cite{Du1}. The automorphism group $\Aut\,V$ of such a surface
$V$ is infinite dimensional and admits an amalgamated free product
structure \cite{DaGi2}.

\subsection{$\C^*$-surfaces with trivial
Makar-Limanov invariant}
Some interesting classes of normal affine surfaces with a trivial
Makar-Limanov invariant were discussed e.g., in \cite{BaML2,
BaML3, DaiRu, Du2} and \cite{MaMiy}. {\it If, for instance, such a
surface $V$ is smooth and its canonical bundle $K_V$ is trivial
(e.g., if $V$ is a smooth complete intersection) then $V\cong
\Spec \C[t,u,v]/(uv-P(t))$ for a polynomial $P\in\C[t]$ with
simple roots} \cite{BaML3} (cf. Example \ref{derplmi}). Here we
concentrate on such surfaces which also admit a $\C^*$-action.
    From Theorems \ref{linact} and
\ref{TN} we deduce:

\bcor\label{epml} A normal affine surface $V$ with an elliptic or
a parabolic $\C^*$-action has a trivial Makar-Limanov invariant if
and only if $V\cong V_{d,e}\cong\A^2_\C/\Z_d$ is an affine toric
surface as in Example \ref{tosu}. \ecor

Actually $V$ as in the corollary admits a parabolic
$\C^*$-action, and so by Remark \ref{adlnd}.1(i) it has a
$\C_+$-action of fiber type and also a $\C_+$-action of
horizontal type (see Examples \ref{tosu} and \ref{pap}).

The following theorem together with Corollary \ref{epml}
describes all normal affine $\C^*$-surfaces with a trivial
Makar-Limanov invariant.

\bthm\label{long} We let $A=A_0[D_+,D_-]$, where $A_0=\C[t]$ and
$D_+$, $D_-$ are $\Q$-divisors on $\A^1_\C$ with $D_++D_-\le 0$.
The following conditions are equivalent.
\begin{enumerate}
\item[(i)] The Makar-Limanov invariant of $V$ is trivial.
\item[(ii)] $A$ admits two homogeneous locally nilpotent derivations
$\p_+,\,\p_-$ of positive and negative degree, respectively, such that
the orbits of the corresponding
$\C_+$-actions are generically different.
\item[(iii)] There are
(not necessarily distinct) points $p_+$, $p_-\in\A^1_\C$ such that
the fractional part $\{D_\pm\}$ of $D_\pm$ is zero or is
supported in $p_\pm$, and $D_++D_-\neq 0$.
\end{enumerate}
\ethm

\bproof The implication $(ii)\Rightarrow (i)$ is evident. For the
proof of the converse, assuming (i) there exist two
non-equivalent locally nilpotent derivations on $A$, which
means that they have different kernels.  By Lemma
\ref{kernel} not both of them can be linear combinations
of derivations of positive degrees, and similarly not both
of them can have homogeneous components of only negative
degree. Thus there are also homogeneous locally nilpotent
derivations on
$A$ of positive and of negative degree. To show
that the corresponding
$\C_+$-actions are not equivalent, we let $v_+$ and $v_-$ be
generators of the $A_0$-modules $A_{d_+}$ and $A_{-d_-}$,
respectively, where $d_+:=d(A_{\ge 0})$ and $d_-:=d(A_{\le 0})$.
By Lemma \ref{kernel} $\ker \p_\pm
=\C[v_\pm,v_\pm^{-1}]\cap A$. Thus, if $\p_+$ and $\p_-$ were
equivalent then $v_\pm$ would be units and so by Corollary
\ref{rd} we would have $A\cong \C[z,v_+,v_+^{-1}]$. As the latter
ring does not admit two non-equivalent $\C_+$-actions, (ii)
follows.

$(iii)\Rightarrow (ii)$. Assuming (iii)  Corollary
\ref{crucol} shows that there are homogeneous derivation $\p_+$
and $\p_-$ of positive and negative degree, respectively. By our
assumption $D_++D_-\ne 0$, hence $A^{\times}=\C$ and so, the
elements $v_+$ and $v_-$ are not units (see \cite[Remark
4.5]{FlZa1}). Thus with the same arguments
as above the derivations
$\p_+$ and $\p_-$ are not equivalent.

$(ii)\Rightarrow (iii)$. Conversely, if (ii)  holds then by
Theorem \ref{main} the first two conditions in (iii) are
satisfied. With the same arguments as above $A$ cannot contain a
non-constant unit, hence again by \cite[Remark 4.5]{FlZa1} we
have $D_++D_-\ne 0$. \eproof

\brem\label{degfib} For explicit equations of $\C^*$-surfaces
with a trivial Makar-Limanov invariant we refer the reader to
Proposition 4.8 in \cite{FlZa1}, where for
$\{-D_\pm\}=e'_\pm/d_\pm[p_\pm]$ one must let
$P_\pm:=(t-p_\pm)^{d'_\pm-e'_\pm}$ with $d'_\pm :=d_\pm/k$ and
$k:=\gcd\,(d_+,d_-)$. \erem

We note that the two locally nilpotent derivations as in Theorem
\ref{long}(ii) do not commute except in the case $V\cong
\A_\C^2$. This is a consequence of the next result. Although it
follows immediately from Lemma \ref{gtl}(a), we provide a direct
argument.

\bcor\label{commut} If a normal affine variety $V=\Spec A$ of
dimension $n$ admits an effective $\C_+^n$-action, then $V\cong
\A_\C^n$. \ecor

\bproof
Let $\A^n_\C\cong \C_+^n.p\hto V$ be an open orbit and consider
the  associated inclusion of $\C$-algebras$A\hto
B:=\C[X_1,\ldots,X_n]$. The derivations
$\p_i:=\frac{\p}{\p X_i}$ on $B$ stabilize $A$ and the
restrictions $\p_i|A$ are the infinitesimal generators of the
actions of the factors of
$\C_+^n$ on $A$. By Proposition \ref{Ren}(b), for every $1\le
i\le n$ the intersection
$$
K_i:=A\cap \bigcap_{j\ne i}\ker\p_i=A\cap \C[X_i]
$$
has transcendence degree 1, hence $K_i\ne \C$. As $\p_i$ acts on
$K_i$ and decreases the degree of polynomials in $K_i$ by 1,
$K_i\subseteq A$ must contain a linear polynomial $a_iX_i+b_i$ and
hence also $X_i$. It follows that $A=B$, as required.
\eproof

For a normal affine surface $V=\Spec A$ with two different affine
rulings $v_+,v_-: V \to \A_\C^1$, Miyanishi and Masuda
\cite{MaMiy} introduced a useful invariant $\iota (v_+,v_-)\in\N$,
called the {\it intertwining number} of $v_+$ and $v_-$, which is
the intersection number of two general fibers of $v_+$ and $v_-$,
respectively. Actually $\iota (v_+,v_-)=\trdeg \left(\Frac A
: \C(v_+,v_-)\right)$.

\bdefi\label{MMI} Let us call the {\it
Miyanishi-Masuda invariant} of $V$ the integer
$$\MM(V):=\min_{(v_+,v_-)} \iota (v_+,v_-)\,,$$
where the minimum is taken over all possible choices of pairs
$(v_+,v_-)$ as above. In case that $V$ is endowed with an
effective $\C^*$-action, we also consider the homogeneous version
$$\MMh(V):=\min_{(v_+,v_-)} \iota (v_+,v_-)\,,$$
under the additional assumption that $v_+$ and $v_-\in A $ as
above are homogeneous. \footnote{Clearly $\MMh\ge\MM (V)$,
where presumably the equality holds.}
\edefi

We let as before $d_+:=d(A_{\ge 0})$ and $d_-:=d(A_{\le 0})$. We
recall \cite[Lemma 3.5]{FlZa1} that $d(A_0[D])$ is equal to the
minimal integer $d\in \Z$ such that the divisor $dD$ is integral.

\blem\label{MMDD} For a
normal affine $\C^*$-surface $V=\Spec A$ with
a trivial Makar-Limanov invariant the following hold.
\begin{enumerate} \item[(a)] If
$A=A_0[D]$ then $\MMh(V)=d(A)$.
\item[(b)] If
$A=A_0[D_+,D_-]$ then $\MMh(V)=-d_+d_-\deg(D_++D_-)$.
\item[(c)] If $\MMh(V)=1$ then $V\cong \A_\C^2$.
\end{enumerate}
\elem

\bproof (a) In this case the grading on $A$ is parabolic, so
    $V$ is a toric surface $V_{d,e}$, where $d=d(A)$, and
the two $\C^*$-equivariant affine rulings on $V$ are provided by
elements $t\in A_0=\C[t]$ and $v\in  A_d=vA_0$ (see
Corollary \ref{epml}). Since the restriction of $v$ onto a general
fiber of $t$ has degree $d$, the result follows.

(b) In this case the grading on $A$ is hyperbolic, and so the two
$\C^*$-equivariant affine rulings on $V$ are provided by elements
$v_\pm\in A_{\pm d_\pm}$ with $A_{\pm d_\pm}=v_\pm A_0$ (see the
proof of Theorem \ref{long}). By Proposition 4.8 in \cite{FlZa1},
$V$ is a cyclic branch covering of degree $k:=\gcd(d_+,d_-)$ of
the normalization of the hypersurface
$\{v_+^{d_-'}v_-^{d_+'}-P(t)=0\}$ in $\A_\C^3=\Spec \C[t,
v_+,v_-]$, where $d'_\pm:=d_\pm/k$. Hence $\MMh(V)=k\deg P(t)$. By
Lemma 4.7 in {\it loc. cit.} we have
$$D_+=D_0+\{D_+\}\quad\mbox{and} \quad D_-=\{D_-\}-D_0-\div\ Q\,,$$
where $Q\in\C[t]$. From (8) and (10) in {\it loc. cit.} we obtain
$$\div\ P= kd_+'d_-'\div\ Q-d_-'\div\ P_+-d_+'\div\
P_-\,,$$ where $\div\ P_\pm=d_\pm\div \{D_\pm\}$. Therefore
\be\label{dip} \div\ P=kd_+'d_-'\div \left(
Q-\{D_+\}-\{D_-\}\right)=-kd_+'d_-'(D_++D_-)\,.\ee Now
$$\MMh(V)=\iota (v_+,v_-)=k\deg
P(t)=-k^2 d_+'d_-'\deg (D_++D_-)=-d_+d_-\deg(D_++D_-)\,,$$ as stated.

(c) The equalities $\MMh(V)=k\deg P(t)=1$ imply
that $k=1$ and $\deg P(t)=1$. Now the assertion follows. \eproof

\subsection{Families of $\C_+$-actions on
a $\C^*$-surface} We show in Corollary \ref{cf} below that any
$\C^*$-surface with a trivial Makar-Limanov invariant admits a
continuous family of generically non-equivalent locally nilpotent
derivations. This is based on the following general observation.

\bprop\label{contfam} If a domain $A$ of finite type admits two
non-commuting locally nilpotent derivations $\p,\,\delta\in \Der
A$, then $A$ also admits a continuous family of generically
non-equivalent locally nilpotent derivations
$\{\p_t\}_{t\in\C^*}\subseteq \Der A$. \eprop

\bproof Letting $\varphi_t=\exp (t\p)$, $\psi_t= \exp
(t\delta)$ be the associated $\C_+$-actions on $A$, we consider
the following two families of conjugated locally nilpotent
derivations on $A$ :
$$\p_t:=\psi_t\circ \p \circ \psi_t^{-1}
\qquad\mbox{and}\qquad \delta_t:=\varphi_t\circ \delta \circ
\varphi_t^{-1}\,.$$ Suppose in contrary that none of these has
the desired property that is, the derivations in each family
$\{\p_t\}_{t\in\A_\C^1}$ and $\{\delta_t\}_{t\in\A_\C^1}$ are
mutually equivalent. It follows that \be\label{me}
\quad\p_t=f(t)\p_0=f(t)\p \quad\mbox{and}\quad
\delta_t=g(t)\delta, \qquad\mbox{where}\qquad f,\,g\in (\Frac
A)[t]\, \ee (see Definition \ref{mli}). Moreover $f(t)\in \ker
\p$ and $g(t)\in \ker \delta$ $\forall t\in\A_\C^1$, where
$f(0)=g(0)=1$. From (\ref{me}) we get:
$$
\psi_t\circ \p=f(t)\p\circ\psi_t
\Rightarrow \frac{\psi_t\circ \p-\p}{t}=
\frac{f(t)\p\circ\psi_t -\p}{t}$$
$$\Rightarrow  \frac{\psi_t-\mbox{id}}{t}\circ
\p=\p\circ
\frac{f(t)\psi_t-\mbox{id}}{t}=\p\circ
\left(
f(t) \frac{\psi_t-\mbox{id}}{t}+\frac{f(t)-1}
{t}\mbox{id}\right)\,.$$
Taking limits as $t\to 0$ we obtain:
$$\delta\circ\p=\p\circ\delta+f'(0)\p\,.$$
Similarly
$$\varphi_t\circ\delta=g(t)\delta\circ\varphi_t\quad\Rightarrow
\quad\p\circ\delta=\delta\circ\p+g'(0)\delta\,.$$
Thus
$$[\p,\,\delta]=a\delta=b\p\,,$$
where $a:=g'(0)\in\ker\delta$ and $b:=-f'(0)\in\ker\p$ (see
Proposition \ref{Ren}(b)). Consequently $\p$ and $\delta$ are
equivalent and hence commute. This contradicts our assumption.
\eproof

\bcor\label{cf} Any normal affine surface $V=\Spec A$ with a
$\C^*$-action and a trivial Makar-Limanov invariant admits
continuous families of $\C_+$-actions and of generically distinct
affine rulings $V\to\A_\C^1$.\ecor

\subsection{Actions
with a big orbit} As an application of our results we give below a
new proof for the classification due to Gizatullin \cite{Gi1} and
Popov \cite{Po}, mentioned in the introduction. Let us recall it
again.

       \bthm\label{pgd} Let a normal affine surface $V$
admits an action
       of an algebraic group $G$ with an open orbit $O$ such that
       $V\backslash O$ is finite.
If $V$ is smooth then $V$ is
       isomorphic to one of the following 5 surfaces:
\be\label{smooth}\A_\C^2,\,\,\,\A^1_\C\times\C^*,
       \,\,\,\C^*\times\C^*,\,\,\,(\pP^1\times\pP^1)\backslash
       \Delta,\,\,\,\pP^2\backslash \bar\Delta\,,\ee where
       $\Delta\subseteq\pP^1\times\pP^1$ is the diagonal
       and $\bar\Delta\subseteq\pP^2$ is a smooth conic.
If $V$ is singular then $V$ is isomorphic to a Veronese cone
$V_{d,1}$ for some $d\ge 2$ (see Example \ref{pop}). \ethm

\brem\label{conic}  Popov \cite{Po} listed as well all affine
surfaces with a big open orbit without the assumption of
normality. \erem

\bproof[Proof of Theorem \ref{pgd}.] We note first that all
surfaces listed in \ref{pgd} admit an action of an algebraic
group with a big open orbit (see Examples \ref{daoa} and
\ref{pop}). Conversely, suppose that $V$ admits an effective
$G$-action with a big open orbit. If $G$ is solvable then by
Lemma \ref{EE}(b) $V$ is isomorphic to $\A^2_\C$, $\A_\C^1\times
\C^*$ or $\C^{*2}$. Otherwise by Lemma \ref{EE}(c)
$G$ contains a subgroup isomorphic to $\SL_2$ or $\PGL_2$. Now
the conclusion follows from the next result.
\eproof

\bprop\label{SL} If $\SL_2$ acts nontrivially on a normal affine
surface $V=\Spec A$ then $V$ is isomorphic either to one of the
surfaces $\pP^1\times \pP^1\backslash \Delta,\quad\pP^2\backslash
\bar\Delta$ or to a Veronese cone $V_{d,1}$. Moreover, any two
such $\SL_2$-actions on $V$ are conjugated in $\Aut (V)$. \eprop

The proof is preceded by the following observations and by Lemma
\ref{pairs} below.

\bsit \label{sl2}
With the assumptions of \ref{SL}, the kernel of
$\SL_2\to\Aut (V)$ is either trivial or equal to the center
$Z(\SL_2)=\{\pm I_2\}$, so one of the groups $G=\SL_2$ or
$G=\PGL_2$ acts effectively on $V$. We let $e=e(G)$ be the
order of the center
$Z(G)$ that is, $e=2$ if $G=\SL_2$ and $e=1$ if $G=\PGL_2$. The
effective $\C^*$-action on $V$ provided by the maximal torus of
diagonal matrices $\T$ of $G$ defines a grading
$A=\bigoplus_{i\in\Z} A_i=A_+\oplus A_0\oplus A_-$. The Borel
subgroups $B_\pm\cong G_e$ (cf. Remark \ref{newder}.2) act
effectively on $V$, and the infinitesimal generators of the
unipotent subgroups $U_\pm\cong\C_+$ of upper/lower triangular
matrices with 1 on the diagonal induce nonzero homogeneous
locally nilpotent derivations $\p_\pm\in\Der\, A$ of degree $\pm
e$ (see Lemma \ref{sdp}). We let $\delta\in\Der\, A$ be the
infinitesimal generator of $\T$ so that $\delta (a)= \deg a\cdot
a$ for $a\in A$ homogeneous. If $\p\in \Der\, A$ is a homogeneous
derivation then $[\delta,\p]=\deg\p\cdot \p$; in particular
$$
[\delta,\p_\pm]=\pm e\p_\pm\quad\mbox{and moreover}\quad
[\p_+,\p_-]=\delta \,.
$$
The adjoint action on $\T$ of the element $\tau=\left({ 0\atop
1}{-1\atop 0}\right)\in G$ of order $2e$  is given by ${\rm
Ad}\,\tau:\delta\mapsto-\delta$. Hence $\tau$ acts on $A$
homogeneously by reversing the grading, i.e.\
$\tau(A_i)=A_{-i}$, and the action of
${\rm Ad}\,\tau$ on the Lie algebra
$\fg\cong \sl_2\cong \C\delta\oplus
\C\p_+\oplus\C\p_-$ of
$G$  is given by $\p_\pm\mapsto
-\p_\mp$. In particular, the $\C^*$-action on $V$
defined by $\T$ is hyperbolic.\esit

\bdefi\label{equipairs} We say that two pairs $(D_+,D_-)$ and
$(\hat D_+, \hat D_-)$ of $\Q$-divisors on $\A_\C^1$ are
equivalent if one can be obtained from the other by applying an
affine transformation $ \A_\C^1\to \A_\C^1$ and a shift
$D_\pm\mapsto D_\pm\pm D_0$ with an integral divisor $D_0$. \edefi

\blem\label{pairs} Let the assumptions be as in Proposition
\ref{SL}. If $A=A_0[D_+,D_-]$ with $A_0=\C[t]$ and $D_++D_-\le 0$
is a DPD representation for $A$ graded via the $\T$-action,
then
$(D_+,D_-)$ is  equivalent to one of the following pairs:
\begin{enumerate}
\item $(0,-[1]-[-1])$, here $e=1$ and $V\cong \pP^1\times
\pP^1\backslash\Delta$, see (\ref{dpd1});
\item $(\frac{1}{2}[0],
  -\frac{1}{2}[0]-[1])$, here $e=1$ and
$V\cong\pP^2\backslash\bar\Delta$, see (\ref{dpd2});
\item $(-\frac{1}{d}[0],-\frac{1}{d}[0])$ with $d\ge 1$; here
$e=1$ and $V\cong V_{2d,1}$, see (\ref{dpd3});
\item $(-\frac{e'}{d}[0], \frac{e'-1}{d}[0])$ with
$d=2e'-1\ge 1$; here $e=2$ and $V\cong V_{d,1}$, see
(\ref{dpd4}).
\end{enumerate}
   \elem

  \bproof
We will first treat the following particular case. {\it If
$D_\pm$ are integral divisors then, in a suitable coordinate
$t$ on $\A_\C^1$, one of the following 3 cases occurs:
\begin{enumerate} \item[($\alpha$)] $e=2$, $D_++D_-=-[0]$.
\item[($\beta$)] $e=1$, $D_++D_-=-2[0]$.
\item[($\gamma$)] $e=1$, $D_++D_-=-[1]-[-1]$.\end{enumerate}
In particular, $(D_+,D_-)$ is equivalent to one of the integral
pairs in (1)-(4).}

\smallskip

To prove this claim, we note first that $D_\pm$ being integral
$A\cong A_0[0,D_++D_-]$ is the normalization of the ring
$B_{1,P}=\C[t,u_+,u_-]/(u_+u_--P(t))$, where $P\in \C[t]$ is a
unitary polynomial with $\div (P)=-(D_++D_-)$ (see Corollary
\ref{KN} and Example \ref{derplmi}). After multiplying $u_+,u_-$
with suitable constants we have
$$\p_\pm (u_\pm)=0,\quad \p_\pm (t)=u_\pm^e,\qquad \mbox{and}\qquad
\p_\pm (u_\mp)=P'(t)P^{e-1}(t)u_\mp^{-e+1}\,$$ (cf. (\ref{dere})).
Hence
$$[\p_+,\p_-](u_+)=\p_+\left(P'(t)P^{e-1}(t)u_+^{-e+1}\right)
=\frac{d}{d t} \left(P'(t)P^{e-1}(t)\right)u_+\,.$$ On the other
hand, $[\p_+,\p_-]=\delta$ and $\delta (u_+)=u_+$, therefore
$$\frac{d}{d t}
\left(P'(t)P^{e-1}(t)\right)=1\,.$$ Thus either $e=2$ and $\deg
P=1$ or $e=1$ and $\deg P=2$. Since $D_++D_-=-\div (P)$, the
special case above follows.

\smallskip

For the rest of the proof we assume that $D_\pm$ are not both
integral. By Theorem \ref{long} (ii)$\Rightarrow$(iii), the
fractional parts
$\{D_\pm\}$ are concentrated in points $p_\pm\in\A_\C^1$.
Clearly, $\tau$ yields an isomorphism
$A_0[D_+,D_-]\cong A_0[\tau_0^*(D_-),\tau_0^*(D_+)]$, where
$\tau_0 : \A_\C^1\to\A_\C^1$ is the affine transformation of
$\Spec A_0=\A_\C^1$ induced by $\tau_0:=\tau|A_0$. By Theorem
4.3(b) in \cite{FlZa1} there is an integral divisor $D_0$ with
\be\label{indi}
  \qquad D_+=\tau_0^*(D_-)+D_0, \quad
D_-=\tau_0^*(D_+)-D_0\,\,\Rightarrow\,\,
D_++D_-=\tau_0^*(D_++D_-)\,.
\ee It follows that
$\tau_0^*(\{D_\pm\})=\{D_\mp\}\neq 0$ and so
$\tau_0(p_\pm)=p_\mp$. With a suitable choice of $t$ then either
(i) $p_+=p_-=0$ is a fixed point of $\tau_0$, or (ii) $\tau:
t\mapsto -t$ and $p_-=-p_+\neq 0$.

We claim that the case (ii) cannot occur. In fact, in
this case we have $\tau_0\ne\id$, and
because of $(\ref{indi})$ and Theorem \ref{main} we may suppose
that
$$
D_+(p_+)=\frac{-e'}{d}, \quad D_-(p_+)=-a,\quad\mbox{and}\quad
D_-(p_-)=\frac{-e'}{d}, \quad D_+(p_-)=-a,
$$
where $d\ge 2$, $0< e'<d$ and $ee'\equiv 1\mod d$. The
points $p_\pm$ are not singular, since otherwise they would
be fixed under the action of the connected group $G$ contradicting
$\tau(p_-)=p_+\ne p_-$. Using \cite[Theorem 4.15]{FlZa1} the
smoothness of the points $p_\pm$ implies that
$e'+ad=1$ and so $e'=1$ and $a=0$. The condition $ee'\equiv 1\mod d$
then forces $e=1$. By Theorem \ref{main}(ii) we also have
$-e(D_+(p_+)+D_-(p_+))\ge 1$, which gives $e(e'/d+a)\ge 1$. This is a
contradiction.

%

Thus in fact (ii) is impossible and so
$p_+=p_-=0$. We can write
$D_+=-\frac{e'}{d}[0]$ and $D_-=-\frac{e'}{d}[0]+E_0$ on
$\A_\C^1=\Spec A_0$ with $d\ge 2$,
where $E_0$ is integral and
$D_++D_-\le 0$. Let $v_\pm\in A_{\pm d}$ be a generator
of $A_{\pm d}$ over $A_0$. Due to Lemmas \ref{cyex},
\ref{DD} and Remark \ref{remDD}
the fraction fields of $A[{\root
d \of v_\pm}]$ and $A[{\root d\of t}]$  are equal,
the normalization $A'$ of $A$ in this field is again
graded, and $\p_\pm$ extend to locally nilpotent
derivations on $A'$. Thus $A'$ admits again an
$\sl_2$-action. Applying Proposition 4.12 in
\cite{FlZa1} $A'\cong A_0'[D_+',D_-']$, where
$A_0'=\C[s]$, $s^d=t$, and $D_\pm'=\pi^*(D_\pm)$ with $\pi
:s\mapsto s^d$. Since the divisors $D_\pm'$ are integral, their
sum
$$D_+'+D_-'=\pi^*(D_++D_-)=-2e'[0]+\pi^*(E_0)$$ is as in
($\alpha$)-($\gamma$) above. In case ($\alpha$) or
($\beta$) clearly $E_0=-b[0]$ with $b\in\Z$. In case
($\alpha$) we have $2e'+db=1$, and since $0<e'<d$ this implies
$b=-1$, so $ d=2e'-1$ and $(D_+,D_-)$ is as in (4).
Similarly, in case ($\beta$) we have
$2e'+db=2$, $e=1$, so $ee'\equiv 1\mod d$ implies
$e'=1$, $b=0$, thus $E_0=0$ and we are in case (3).

  In the remaining case ($\gamma$) we have $e=1$.
Letting $E_0=-b[0]+E_0'$ with
  $E_0'(0)=0$, we obtain that
  $-(2e'+db)[0]+\pi^*(E_0')=-[p]-[q]$ with $p\neq q$. Therefore either
  $$\qquad 2e'+db=0\quad\mbox{and}\quad \pi^*(E_0')=-[p]-[q]
  \quad\mbox{with}\quad p,q\neq 0\,,$$ or, up to interchanging $p$ and
  $q$,
  $$\qquad p=0,\quad 2e'+db=1\quad\mbox{and}\quad
  \pi^*(E_0')=-[q]\neq [0]\,.$$ Actually this latter case cannot
  occur since $d\ge 2$ divides $\deg \pi^*(E_0')$. Thus we must
  have $d=2$, $e'=1$, $b=-1$ and $p=-q\neq 0$. Letting e.g., $p=1$
  we obtain that $(D_+,D_-)$ is as in (2). This proves the
lemma.
\eproof

\bproof[Proof of Proposition \ref{SL}.]
Lemma \ref{pairs} implies that a surface with an
$\SL_2$-action is isomorphic to one of the surfaces listed in
the proposition. It remains to show that this isomorphism can
be chosen to be equivariant with respect to
the given $\SL_2$-actions. For this we restrict to the case
$\pP^2\backslash \bar\Delta$, the argument in the
other cases being similar.

Let $V=\Spec A$ be an $\SL_2$-surface as in Lemma
\ref{pairs}(1) and denote $V':=\pP^2\backslash\bar\Delta$
with its standard action as in Example \ref{daoa}. Both $A$
and the affine coordinate ring
$A'$ of $V'$ are equipped with the grading
coming from the maximal torus in $\SL_2$, and by the
construction in  Lemma \ref{pairs} the isomorphism $A\cong
A'$ is compatible with these gradings. Let
$(\delta,\p_+,\p_-)$ be the triplet of
derivations on
$A$ as in \ref{sl2}, and let $(\delta',\p'_+,\p'_-)$ denote the
corresponding derivations on $A'$. Using
Lemma \ref{pairs} again $e=\deg\p_\pm=\pm 1$; as
$\PGL_2$ acts on $V'$ (cf.\ Example \ref{daoa}) we
also have $\deg \p_\pm'=\pm1$.

Now Proposition \ref{unique} shows that the pairs
$(A,\p_+)$ and
$A',
\p'_+)$ are isomorphic, so there is a graded
isomorphism $f:A\to A'$ with $f_*(\delta)=\delta'$ and
$f_*(\p_+)=\p_+'$. Again by Proposition \ref{unique}
$f_*(\p_-)=c\p'_-$ for some constant $c\in\C^*$. As
$\delta=[\p_+,\p_-]$ it follows that
$\delta'=f_*(\delta)=f_*([\p_+,\p_-])=c[\p_+',\p_-']=c\delta'$.
Hence
$c=1$ and so $f_*(\p_\pm)=\p_\pm'$.  By Proposition
\ref{groupder} this means that the induced  isomorphism
$V\cong V'$ is equivariant with respect to the Borel
subgroups $B_\pm$ of $\SL_2$ and so it is
$\SL_2$-equivariant, as desired.
\eproof

\brem\label{linsl2} Proposition \ref{SL} shows in particular
that any $\SL_2$-action on the plane $\A_\C^2$ is conjugated in
$\Aut\ \A_\C^2$ to the standard linear representation.\erem

\section{Concluding remarks. Examples}

Here we illustrate our methods in concrete examples. According to
Gizatullin's Theorem cited in \ref{pgd}, there are only 5
different homogeneous affine surfaces (\ref{smooth}). In the
following example we consider more closely the last two of these
surfaces $\pP^1\times\pP^1\backslash
       \Delta$ and $\pP^2\backslash
       \bar\Delta$  (cf. \cite[Lemma 2]{Po}).

\bexa\label{daoa}
Let $V\cong \C^2$ be a 2-dimensional vector space. The group
$\PGL_2\cong
\PGL(V)$ then acts on $\pP^1=\pP(V)$ as well as on the
projectivized space of binary quadrics
$\pP^2=\pP(S^2V)$. Since
$\PGL_2$ acts doubly transitive on $\pP^1$, the diagonal action
on
$\pP^1\times\pP^1$ has an open orbit
$\pP^1\times\pP^1\backslash\Delta$, where $\Delta$ is the diagonal.
Similarly, the action of $\PGL_2$ on $\pP^2$ leaves the
degenerate quadrics invariant thus providing an action on
$\pP^2\backslash
\bar\Delta$, where $\bar\Delta$ is the space of degenerate binary
forms.

The symmetric product $V\times V\to S^2V$, $(v,w)\mapsto v\vee w$,
induces a natural unramified 2:1 covering
$$
p:\pP^1\times \pP^1\backslash \Delta\lto \pP^2\backslash \bar\Delta,
$$
where the covering involution is the map interchanging the two
factors of   $\pP^1\times \pP^1$.

To make the situation more explicit, fix a basis $v_0,v_1$ of
$V$ so that the points of $\pP(V)$ can be represented in
coordinates $[x_0,x_1]$. With respect to the basis $v_0^2$,
$2v_0v_1$, $v_1^2$ of $S^2V$ the points of $\pP^2=\pP(S^2V)$
have then  coordinates
$[u'_+,s,u'_-]$. Clearly $\bar\Delta$ has equation
$Q:=s^2-u'_+u'_-=0$.  The map
$p$ factors through
\begin{diagram}
\pP^1\times \pP^1\backslash \Delta&\rTo^{\tilde p}&
H&\rTo^{can}&
\pP^2\backslash \bar\Delta,
\end{diagram}
where $H$ is the affine quadric $\{Q=1\}\subseteq \C^3$ and
$\tilde p$ is the isomorphism given by
$$
([x_0,x_1], [y_0,y_1])\mapsto \frac{1}{x_0y_1-x_1y_0}
(2x_0y_0, x_0y_1+x_1y_0, 2x_1y_1).
$$
This isomorphism identifies the factors interchanging involution of
$\pP^1\times \pP^1$ with the map $(u_+,s,u_-)\mapsto
-(u_+,s,u_-)$.  Thus $\pP^1\times \pP^1\backslash
\Delta\cong H\cong \Spec A'$, where according to Example
4.10 in \cite{FlZa1}
\be\label{dpd1}
A':=\C[u'_+,s,u'_-]/(u'_+u'_--s^2+1)\cong A'_0[D'_+,D'_-]
\ee
with $A'_0:=\C[s]$, $D'_+=0$ and $D'_-=-[1]-[-1]$. This isomorphism
determines a hyperbolic grading with $\deg s=0$ and $\deg u'_\pm=\pm
1$.

Next we turn to the surface $\pP^2\backslash\bar\Delta$, which
is the spectrum of the invariant ring $A:= A^{\prime\Z_2}$.
As noted above the action of $\Z_2$ on $A'$ is given by
$(u_+',s,u_-')\mapsto -(u_+',s,u_-')$. The algebra of
invariants is generated by the degree 2
monomials in
$s,u'_\pm$
$$
u_\pm:=su'_\pm,\quad v_\pm:=u_\pm^{\prime
2}\quad\mbox{and}\quad
   t:=s^2
$$
satisfying the relations
$$
u_-=t(t-1)u_+^{-1},\quad
v_+=t^{-1}u_+^{ 2}, \quad v_-=t(t-1)^2 u_+^{-2}\,
$$
(observe that
$u_+'u_-'=s^2-1=t-1$ in $A'$). Thus $A=\C[t][v_-,u_-,u_+,v_+]$
can be presented as
\be
\label{sva}\qquad A=\C[t]\left[t(t-1)^2 u_+^{-2},
\,t(t-1)u_+^{-1},\,u_+,\,t^{-1}u_+^2\right]\subseteq
\C(t)[u_+,u_+^{-1}]\,.
\ee
By virtue of (\ref{sva}) and Lemma 4.6 in
\cite{FlZa1},
\be\label{dpd2}
A\cong A_0[D_+,D_-] \quad \mbox{with}\quad
D_+=\frac{1}{2}[0],\quad D_-=-\frac{1}{2}[0]-[1];
\ee
indeed, according to this lemma
$$
D_+=-\min \left\{0,-\frac{1}{2}[0]\right\}=\frac{1}{2}[0],
\,\,\, D_-=
-\min \left\{\div\ t(t-1),\frac{\div\ t(t-1)^2}{2}
\right\}=-\frac{1}{2}[0]-[1]
$$
and so $D_++D_-=-[1]$.

With this example one can also make some of the previous results quite
explicit. For instance, $\pi^*(D_+)=D_+'+D_0$ with $D_0:=[0]$ and
$\pi^*(D_-)=-[0]-[1]-[-1]=D_-'-D_0$ with
$\pi: s\longmapsto s^2=t$, which agrees with Proposition 4.12 in
\cite{FlZa1} applied to the Galois $\Z_2$-extension $A\hookrightarrow
A'$. Further, the fractional parts $\{D_\pm\}$ of $D_\pm$ are
supported at one point; compare with Theorem
\ref{long} above.

For every $\lambda=(\lambda_+, \lambda_0,\lambda_-)$ with
$\lambda_0^2=4\lambda_+\lambda_-$ the hyperplane in $\pP^2$ given by
$f_\lambda:=\lambda_0s+\lambda_+u_++\lambda_-u_-=0$ intersects
$\bar\Delta$ in one point. It follows that the maps
$$
f_\lambda: H\lto \A^1_\C\quad\mbox{and}\quad
g_\lambda :=f^2_\lambda/Q: \pP^2\backslash \bar\Delta\lto \A^1_\C
$$
provide explicit families of affine rulings compatible with $p$
(cf.\  Proposition \ref{contfam}). By \cite[Proposition
1.11]{Be2} any affine ruling $\pP^2\backslash \bar\Delta\to
\A^1_\C$ is given by a certain $g_\lambda$; they can be
visualized via the Segre and Veronese embeddings
$\pP^1\times\pP^1\hto \pP^3$, $\pP^2\hto \pP^5$.

Finally it is easy to see (and left as an exercise to the reader)
that the locally nilpotent derivations $\p_\pm$ defined by the
unipotent subgroups $U_\pm\subseteq \PGL_2$ of upper/lower
triangular matrices with 1 on the diagonal are of degree 1 and
are given by the formulas in Remark
\ref{mulpar}(1) (compare also with the proof of Lemma
\ref{SL}).

\eexa

\bexa\label{pop} {\it Veronese cones.} For $d\ge 1$ and $e=1$,
$A_{d,1}=\bigoplus_{\nu\ge 0} \C[X,Y]_{\nu d}$ is the $d$-th
Veronese subring of the polynomial ring $\C[X,Y]$. The standard
$\SL_2$-action on $\C[X,Y]$ stabilizes $A_{d,1}$ and so, induces
an $\SL_2$-action on the normal affine surface
$V_{d,1}:=\Spec\,A_{d,1}$. This $\SL_2$-action has a unique fixed
point $\bar 0\in V_{d,1}$ and is transitive on $V_{d,1}\backslash
\{\bar 0\}$.

The algebra $A_{d,1}$ is generated by the monomials $X^iY^{d-i}\in
(A_{d,1})_1$ $(i=0,\ldots,d)$, and these define an embedding
$\rho : V_{d,1}\hookrightarrow \A_\C^{d+1}$ onto the affine cone
over the degree $d$ rational normal curve $\Gamma_d\cong \Proj
A_{d,1}$ in $\pP^d$. The morphism $\rho$ is equivariant w.r.t.
the standard irreducible representation of $\SL_2$ on the space
$\A_\C^{d+1}$ of degree $d$ binary forms. The group $\SL_2$
(respectively, $\PGL_2$) acts  effectively on $V_{d,1}$ if $d$ is
odd (respectively, even). The stabilizer subgroup
$$
N_d:=\Bigl\{\left(\begin{array} {cc}
\varepsilon & \alpha\\
0 & \varepsilon^{-1}
\end{array}\right)\Bigl\vert\,
\varepsilon^d=1,\,\,\alpha\in\C_+\Bigr\}\,
$$ of the binary form $X^d\in A_{d,1}$ is a cyclic extension of the maximal
connected unipotent subgroup $N=N_1$ of $\SL_2$. Clearly,
$V_{d,1}\cong \Spec\,\cO\left(\SL_2/N_d\right)$, as
$\SL_2/N_d\cong V_{d,1}\backslash \{\bar 0\}$ and $V_{d,1}$ is
normal \cite{Po}.

To represent the Veronese cones via the DPD construction, note
first that the action of the torus $\T= \left\{\left({ \lambda
\atop 0}{0\atop \lambda^{-1}}\right)\ |\
\lambda\in\C^*\right\}\subseteq\SL_2$ provides a grading on the
ring $A=\C[X,Y]$ with $\deg X=1,\, \deg Y=-1$, and so induces a
grading on the $d$-th Veronese subring
$A_{d,1}=A^{(d)}=\bigoplus_{i\in\Z} A^{(d)}_i$. We consider
separately the cases that $d$ is even or odd.

\smallskip

(1) For $d=2d'$ even, the $\T$-action on $A^{(d)}$ factorizes
through an action of $\T/\Z_2$, which corresponds to letting $\deg
X=1/2,\, \deg Y=-1/2$. With $t:=(XY)^{d'}\in A^{(d)}_0=\C[t]$ and
$u:=X/Y\in (\Frac A_0) A^{(d)}_1$, we have
$$u_i:=X^{d'+i}Y^{d'-i}=tu^i\in A^{(d)}_i,\quad -d'\le i \le
d'\,.$$ As $A^{(d)}=\C[u_{-d'},\ldots,u_{d'}]$ by Lemma 4.6
in \cite{FlZa1} $A^{(d)}\cong A_0^{(d)}[D_+,D_-]$, where
\be\label{dpd3}
\quad D_+=-\min_{1\le i \le d'}
\left\{-\frac{1}{i}[0]\right\}=-\frac{1}{d'}[0],\quad D_-=
-\min_{-d'\le i \le -1}
\left\{-\frac{1}{-i}[0]\right\}=-\frac{1}{d'}[0]
\ee
and so $D_++D_-=-\frac{2}{d'}[0]$.

\smallskip

(2) For $d=2e'-1$ odd, the torus $\T$ acts effectively on
$A^{(d)}$. We let $t:=(XY)^d\in A^{(d)}_0=\C[t]$,
$u_1:=X^{e'}Y^{e'-1}\in A^{(d)}_1$ and
$$u_{2k-1}:=X^{e'+k-1}Y^{e'-k}=t^{-k+1}u_1^{2k-1}\in
A_{2k-1}^{(d)},\quad -e'+1\le k\le e'\,.$$ As
$A^{(d)}=\C[u_{-d},\ldots,u_{d}]$ then by Lemma 4.6 in
\cite{FlZa1} $A^{(d)}\cong A_0^{(d)}[D_+,D_-]$, where
$$D_+=-\min_{1\le k \le e'}
\left\{-\frac{k-1}{2k-1}[0]\right\}=\frac{e'-1}{d}[0],\quad D_-=
-\min_{-e'+1\le k \le 0}
\left\{\frac{k-1}{2k-1}[0]\right\}=-\frac{e'}{d}[0]$$ and so,
$D_++D_-=-\frac{1}{d}[0]$. We notice that
\be\label{dpd4}
(D_+,D_-)=\left(\frac{e'-1}{d}[0],-\frac{e'}{d}[0]\right)
\sim \left(-\frac{e'}{d}[0],\frac{e'-1}{d}[0]\right)
\ee
via the shift $(D_+,D_-)\mapsto (D_+-[0],D_-+[0])$.

Alternatively, the Veronese cone $V_{d,1}$ can be obtained from
the Hirzebruch surface $\Sigma_d:=\pP \left(\cO_{\pP^1}\oplus
\cO_{\pP^1}(-d)\right)$ by deleting a section $C_d$ with
$C_d^2=d$ and contracting the exceptional section $E_d$ with
$E_d^2=-d$ \cite[\S 11, Example 1]{DaGi2}. This leads
\cite{DaGi2} to a description of the automorphism groups $\Aut
(V_{d,1})$.
\eexa

In the next example we exhibit affine surfaces $V$ such that the
automorphism group $\Aut\ V$ acts on $V$ with a big open orbit
$O$ and there are algebraic group actions on $V$ with an open
orbit, whereas there is no such action with a big open orbit.

\bexa\label{newgr} ({\it Actions on surfaces with a big open
orbit.}) Let $D_\pm$ be two $\Q$-divisors on $\A^1_\C$ with
$D_++D_-\le 0$ such that the supports of the fractional parts
$\{D_\pm\}$ are contained in (possibly the same) points
$\{p_\pm\}$. According to Theorem \ref{long} the ring
$A:=A_0[D_+,D_-]$ with $A_0:=\C[t]$ admits locally nilpotent
derivations $\p_\pm$ of positive and negative degree. The
associated $\C_+$-actions $\varphi_+$ and $\varphi_-$ on $V$ are
not equivalent provided that $D_++D_-\ne 0$ (see Definition
\ref{mli}).

Consider the subgroup $G:=\langle
\varphi_+,\lambda,\varphi_-\rangle \subseteq\Aut\,V$ generated by
$\varphi_\pm$ and the $\C^*$-action $\lambda$ on $V$. The fixed
points set of $G$ is finite as it is contained in the fixed points
set $F$ of the $\C^*$-action on $V$. Recall that $F$ has exactly
one point $a'$ over every point $a\in \A_\C^1$ with
$D_+(a)+D_-(a)<0$ \cite[Theorem 4.18(b)]{FlZa1}. We claim that $G$
acts transitively on the complement $V\backslash F$. Indeed, the
algebraic subgroup $G_{e_\pm}:=\langle
\varphi_\pm,\lambda\rangle$ of $G$ acts on $V$ with an open orbit
which contains $V\backslash v_\pm^{-1}(0)$. Hence for a general
point $x\in V$, the orbit $G . x$ contains
$V\backslash\{v_+^{-1}(0)\cap v_-^{-1}(0)\}=V\backslash F$ (cf.
Proposition \ref{qmap}).

Thus $G$ acts on $V$ with a big open orbit. However, such a
surface $V$ does not admit an action of an {\em algebraic} group
with a big open orbit unless it is isomorphic to one of the
surfaces from Theorem \ref{pgd}. For instance, this is the case
if $V$ has two or more singular points (cf. \cite[Theorem
4.15]{FlZa1}), or is an affine toric surface $V_{d,e}$ with $d>e
>  1$.

A particular case is provided by the dihedral surfaces
$V_{d,d-1}\cong \Spec A_{d,d-1}$, where $A_{d,d-1}\cong
\C[t,u_+,u_-]/(u_+u_--t^d)$ and $d\ge 3$.  We have
$A_{d,d-1}\cong_{A_0} A_0[D_+,D_-]$ with $D_+=0$ and $D_-=-d[0]$
for a grading on $A_{d,d-1}$ with $\deg t=0,\,\,\deg u_\pm=\pm 1$
(see Corollary \ref{KN}). The derivations
$$
\p_\pm=u_\pm\,\p/\p t+d\cdot t^{d-1}\p/\p u_\mp
$$
with $\deg \p_\pm=\pm 1$ are locally nilpotent on $A_{d,d-1}$. The
associated $\C_+$-actions $\varphi_+$ and $\varphi_-$ on
$V_{d,d-1}$ generate a subgroup $G$ of $\Aut \,V_{d,d-1}$. Using
e.g.,\ Remark \ref{aut}.4 it is easily seen that $G$ acts with a
big open orbit $V_{d,d-1}\backslash\{\bar 0\}$, where $\bar 0\in
V_{d,d-1}$ denotes the unique singular point. The dihedral
surfaces $V_{d,d-1}$ with $d\ge 3$ are not isomorphic to Veronese
cones, since the exceptional set of the minimal resolution of
$V_{d,d-1}$ is a chain of $d-1$ rational $(-2)$-curves, whereas
it is just one rational curve for every singular Veronese cone.
Hence by Popov's Theorem  \ref{pgd} there is no algebraic group
action with a big open orbit on $V_{d,d-1}$.
   \eexa

We continue with examples that illustrate Corollaries \ref{KN}
and \ref{maincor}.

\bexa\label{Dani} {\it Danielewski's surfaces.}
These are the
smooth surfaces
$$
W_d:=\{u^dv=t^2+t\}\subseteq\A^3_\C
\qquad (d\ge 1)\,.
$$
Thus $W_d=\Spec\,B_{d,P}$ with $P(t):=t^2+t$ is one of the
surfaces studied in Corollary \ref{KN}. So it admits a
$\C^*$-equivariant $\C_+$-action along the fibers of the affine
ruling $u : W_d\to \A^1_\C$. Note that $W_1\cong
(\pP^1\times\pP^1)\backslash \Delta$ has a continuous family of
affine rulings over $\A_\C^1$ (see Example \ref{daoa}), whereas
for every $d\ge 2$, such a ruling on $W_d$ is unique and
ML$(W_d)=\C[u]$. The latter follows from Theorem \ref{long} as
$A=B_{d,P}\cong A_0\left[D_+,D_-\right]$ with $A_0= \C[t]$,
$D_+=0$ and $D_-=-\frac{1}{d} \left([0]+[-1]\right)$, where the
fractional part $\{D_-\}$ of the $\Q$-divisor $D_-$ is supported
at two points (see Example 4.10 in \cite{FlZa1}).

According to Corollaries 4.24 and 4.25 in \cite{FlZa1} we have
Pic$(W_d)\cong \Z$ generated e.g., by $[\bO_0^-]$, whereas
$K_{W_d}=0$.

We recall \cite{Dani} \footnote{Cf. also \cite{BaML2, Fi, Wi}.}
that these surfaces provide examples of non-cancellation, that is
$W_d\times\A^1_\C\cong W_{d'} \times\A^1_\C\quad \forall d,d'\in
\N$, whereas $W_d\not\cong W_{d'}$ if $d\neq d'$. \eexa

\bexa\label{Bertin} {\it Bertin's surfaces.} These are the smooth
affine surfaces \be\label{Be}
W_{d,n}:=\{x^dy=x+z^n\}\subseteq\A_\C^3\,; \ee they admit an
algebraic group action with an open orbit \cite{Be2}. Note that
$W_{d,1}\cong \A_\C^2$ and $W_{1,n}\cong V_{n,n-1}$ admit
continuous families of affine rulings over $\A_\C^1$. Thus we will
suppose in the sequel that $d,n\ge 2$. The defining equation of
$W_{d,n}$ is quasihomogeneous with weights
$$
\deg x=n, \quad \deg y=-n(d-1), \quad \deg z=1.
$$
To compute a DPD presentation of the coordinate ring
$A=\C[x,y,z]/(x^dy-x-z^n)$, we note that $A_0=\C[t]$ with
$t:=x^{d-1}y-1$. Moreover the equations $xt=z^n$ and
$y=(t+1)t^{d-1}z^{-n(d-1)}$ show that $A=A_0[z,
t^{-1}z^n,(t+1)t^{d-1}z^{-n(d-1)}]$, and so by \cite[Lemma
4.6]{FlZa1}
$$
A\cong A_0[D_+,D_-]\quad\mbox{ with }\quad
D_+=\frac{1}{n}[0]\quad\mbox{ and }\quad
D_-=-\frac{1}{n}[0]-\frac{1}{n(d-1)}[-1]\,.
$$
A homogeneous locally nilpotent derivation $\p$ on $A$ of degree
$nd-1$ can be given by
$$
\p(x)=0,\quad \p(y)=nz^{n-1},\quad \p(z)=x^d.
$$
According to Corollary \ref{maincor} Bertin's surfaces can be
described as cyclic quotients of $V_{d',P}$ for a suitable pair
$(d',P)$. To find such a presentation one takes the normalization
$A'$ of $A$ in the quotient field of $A[u]$, where $u:=x^{1/n}$.
The equation $z^n/u^n=t$ shows that $s:=z/u\in A'$. Thus $A'$
contains $\C[s,u,y]/(u^{n(d-1)}y-1-s^n)$, and since the latter
ring is normal, these two rings are equal. The derivation $\p$
extends to $A'$ via $\p(u)=0$ and $\p(s)= u^{nd-1}$ commuting
with the homogeneous $\Z_n$-action on $A'$
$$
\zeta . s= \zeta^{-1}s,\qquad \zeta . u= \zeta u
\qquad\text{and}\qquad \zeta . y= y\,,
$$
  where $\zeta$ is a primitive $n$-th root of unity. This
action on $V'=\Spec A'$ is fixed point free and
$A=A^{\prime\Z_n}$ i.e., $W_{d,n}\cong V'/\Z_n=V_{d',P}/\Z_n$,
where $d':=n(d-1)$ and $P:=s^n+1$.

For every
$d,n\ge 2$ the fractional part $\{D_-\}$ of the $\Q$-divisor $D_-$
is supported at two points. Hence according to Theorem
\ref{long}, $x: W_{d,n}\to \A_\C^1$ gives a unique affine ruling
on $W_{d,n}$ over an affine base, and ML$(W_{d,n})=\C[x]$ (cf.
\cite{ML2}). The latter also follows from
    \cite[Theorem 1.8 and Example 4.11(iii)]{Be2}
(cf. Theorem \ref{zigzag})
    due to the
fact that the dual graph of a minimal compactification of
$W_{d,n}$ is not linear.

It can be readily seen that Pic$\,(W_{d,n})\cong \Z/n\Z$ generated
e.g., by $[O_0]$, whereas $K_{W_{d,n}}=0$ (see e.g., Corollaries
4.24 and 4.25 in \cite{FlZa1}). \eexa

\brem\label{last} Any affine surface $V\not\cong \A^2_\C$ which
admits an elliptic $\C^*$-action is singular. If $V$ is equipped
with a parabolic $\C^*$-action and a horizontal $\C_+$-action
then by Theorem \ref{TNC} it has a quotient singularity. Thus
being smooth the surfaces $\pP^1\times\pP^1\backslash\Delta$,
$\pP^2\backslash\bar\Delta$, $W_d$ and $W_{d,n}$ with $d,n\ge 2$
admit neither elliptic nor parabolic $\C^*$-actions.

\erem

\end{document}